\input amstex
\documentstyle{amsppt}
\NoBlackBoxes

\loadbold
\magnification = \magstep 1

\define\Dom{\operatorname{Dom}}
\define\Id{\operatorname{Id}}
\define\supp{\operatorname{supp}}
\define\rank{\operatorname{rank}}
\define\hlb{holomorphic line bundle }

\define\tE{\widetilde E}
\define\tX{\widetilde X}
\define\tf{\widetilde F}\define\tp{\widetilde P}\define\tq{\widetilde Q}
\define\tO{\widetilde\Omega}\define\Om{\Omega}
\define\tg{\widetilde G}\define\tm{\widetilde M}				
\define\tM{\widetilde M}\define\te{\widetilde E}\define\p0{P_0}

\define\g{\Gamma}
\define\gn{\gamma_\nu}
\define\gpn{\gamma^\prime_\nu}
\define\gsn{\gamma^{\prime\prime}_\nu}

\define\De{\Cal D}
\define\f{\varphi}\define\lam{\lambda}
\define\db{\bar\partial}
\define\pa{\partial} 
\define\rg{\operatorname{Reg}\,(X)}
\define\sg{\operatorname{Sing}\,(X)}
\define\tl#1{\widetilde #1}
\define\sn#1{\operatorname{Sing}(#1)}
\define\rl#1{\operatorname{Reg}(#1)}
\define\sG{\scriptscriptstyle\Gamma}
\define\ve{\varepsilon}

\define\lat{\widetilde\Delta^{\prime\prime}}
\define\latk{\widetilde\Delta^{\prime\prime}_{k}}
\define\latd{\widetilde\Delta^{\prime\prime}_{k}\upharpoonright_{\tO}}
\define\ng#1{N^k_{\scriptscriptstyle\Gamma}(\lambda,#1)}
\define\dg#1{\dim_{\scriptscriptstyle\Gamma}#1}
\define\nlp{N_{\scriptscriptstyle\Gamma}(\lambda,\widetilde P)}
\define\np0{N(\lambda,P_0)}
\define\n0#1{N^k_U(\lambda,#1)}
\define\nom#1{N^k_{\Omega}(\lambda,#1)}

\topmatter
\title
$L^2$--Riemann--Roch inequalities for covering manifolds
\endtitle

\author{Radu Todor, Ionu\c t Chiose, George Marinescu}
\endauthor

\address{Faculty of Mathematics, University of Bucharest, Str Academiei 14, Bucharest, Romania}
\endaddress
\address{Department of Mathematics, SUNY at Stony Brook, Stony Brook, NY 11794-3651 USA}
\endaddress
\email{chiose\@ math.sunysb.edu}\endemail
\address{Institute of Mathematics of the Romanian Academy, PO Box
1--764, RO--70700, Bucharest, Romania } \endaddress
\curraddr Institut f\"ur Mathematik, Humboldt--Universit\"at zu Berlin,
 Unter den Linden 6,
10099 Berlin, Deutschland 
\endcurraddr

\email{george\@ mathematik.hu-berlin.de}\endemail

\thanks{The main part of this work was carried out while the third named author was supported 
by a DFG Stipendium at the ``\,Graduiertenkolleg
Geometrie und Nichtlineare Analysis\,'', Humboldt--Universit\"at zu Berlin.}\endthanks
\keywords{von Neuman dimension, Moishezon manifold, semi--classical estimate,
Shubin's IMS localization, Demailly's asymptotic formula, Nori's weak Lefschetz theorem}
\endkeywords
\subjclass{32J20}\endsubjclass
\abstract
We study the existence of $L^2$ holomorphic sections of invariant line 
bundles over Galois coverings of Zariski open sets in Moishezon manilolds.
We show that the von Neuman dimension of the space of $L^2$ holomorphic sections
is bounded below under reasonable curvature conditions.
We also give criteria for a a compact complex space with isolated singularities
and some related strongly pseudoconcave manifolds to be Moishezon. Their coverings 
are then studied with the same methods.
As applications we give weak Lefschetz theorems using the Napier--Ramachandran
proof of the Nori theorem.
\endabstract
\toc
\subhead
\S 1 Estimates of the spectrum distribution function
\endsubhead
\subhead
\S 2 Geometric situations
\endsubhead
\subhead
\S 3  Coverings of some strongly pseudoconcave manifolds 
\endsubhead
\subhead
\S 4 $L^2$ generalization of a theorem of Takayama 
\endsubhead
\subhead
\S 5 Further remarks
\endsubhead
\subhead
\S 6 Weak Lefschetz theorems
\endsubhead

\endtoc
\endtopmatter

\document

In this paper we wish to address the following problem. Let $\tM$ be a complex manifold and
assume there is a discrete group $\g\subset\operatorname{Aut}\tM$ acting freely and properly
discontinuously on $\tM$. Suppose that the quotient $M=\tM/\g$ is a Zariski open set in a 
Moishezon manifold $X$ and let $E\longrightarrow X$ be a \hlb on $X$. We denote by $p:\tM
\longrightarrow M$ the canonical projection. 

\definition{Problem}
Find non-trivial $L^2$ holomorphic sections in $p^*E^k$ over $\tM$ for large $k$ provided
$E$ satisfies reasonable conditions in terms of curvature positivity.
\enddefinition

Let us describe briefly the background of this problem. In two earth-breaking papers
Siu \cite{Si1, Si2} proved the Grauert--Riemenschneider conjecture \cite{GR} by showing that
if $E\longrightarrow X$ is a semipositive \hlb on a compact manifold and it is positive at 
one point then $E^k$ has a lot of holomorphic sections i.e. we have the Riemann--Roch inequality
$\dim H^0(X,E^k)\geqslant C\,k^n$
for large $k$, where $n=\dim X$. Demailly \cite{De1} developped a more 
powerful method based on Witten's work \cite{Wi} to get asymptotic Morse inequalities.
Takayama \cite{Ta} generalized the Riemann--Roch inequality for the case 
when $X$ is a compact complex space and $E$ is positive in the neighbourhood of 
an analytic subset. In order to generalize these results to the case of coverings we 
shall use the framework of Atiyah \cite{At} who computed the von Neumann index of 
a $\g$--invariant elliptic operator. Our main technical device comes from a paper of 
Shubin \cite{Sh} in which a proof in the spirit of Witten of the Novikov--Shubin 
inequalities is given.  The present paper 
pertains also to the work of Gromov, Henkin and Shubin \cite{GHS} in which the authors compute
the von Neumann dimension of the space of $L^2$ holomorphic functions on coverings of strictly 
pseudoconvex domains. The von Neumann dimension turns out to be infinite 
generalizing thus Grauert's theorem.

Let us describe the content of our paper. In \S 1 we generalize the Weyl type formula of 
Demailly by describing the asymptotic behaviour 
of the spectrum of a $\g$--invariant laplacian associated to high powers of a 
$\g$--invariant line bundle. Using this tool we prove in \S 2 the main theorem which consists of studying
manifolds with pointwise bounded torsion admitting an uniformly positive line bundle outside a compact set.
In \S 3 we consider a special case of 1--concave manifolds and stongly pseudoconcave domains
associated to compact complex spaces with isolated singularities. Our results are meaningful even for 
the trivial covering and they extend the Demailly--Siu criteria for algebraicity from the case 
of compact manifolds. For the type of manifolds under discussion we also prove some stability results
for the perturbation of complex structure.
In \S 4 we generalize the $L^2$ Riemann--Roch inequality of Takayama \cite{Ta} by considering
Galois coverings of smooth Zariski open sets in compact Moishezon spaces. 
We also remark that by using the 
$\db$--method as in Napier and Ramachandran \cite{NR} we may extend the result for arbitrary coverings.
Paragraph \S 5 gives applications to the quotients of bounded domains in $\Bbb C^n$. 
We remark that the von Neuman dimension of the space of $L^2$ holomorphic pluricanonical
sections is infinite if the volume of the quotient in the Bergman metric is infinite.
At the opposite side we give a positive partial answer to a question of Griffiths, by showing
that the Bergman volume of the quotient is finite provided the quotient is the regular part of 
a compact complex space with only isolated singularities.
Finally, \S 6 is devoted to proving weak Lefschetz theorems for Moishezon manifolds using 
the analytic proof (and generalization) of Nori's results due to \cite{NR}.
\definition{Aknowledgements}
We want to express our wholehearted thanks the following people:
V. Iftimie for iniatiating this project, J. Koll\'ar for bringing the work of 
Takayama to our attention and G. Henkin for useful conversations. The third named author expresses
its gratitude to Prof. J. Leiterer for excellent working conditions and the 
``\,Graduiertenkolleg Geometrie und Nichtlineare Analysis\,'', especially Prof. Th. Friedrich,
for support.
\enddefinition

\subhead
{\S 1. Estimates of the spectrum distribution function}
\endsubhead
Let $\tm$ be a complex analytic manifold of complex dimension $n$ on which a
discrete group
$\Gamma$ acts freely and properly discontinuously. Let $X=\tm/\Gamma$ 
let $\pi:\tm\longrightarrow X$ be the canonical projection. 
We assume $\tm$ paracompact so that $\g$ will be countable. 
Suppose we are given a holomorphic vector bundle $F$ on $X$ and take its pull-back $\tf=\pi^{\ast}F$, 
which is a $\g$ invariant bundle on $\tm$.
We also fix a $\g$ invariant hermitian metric on $\tm$ and on $\widetilde F$. 
We consider a relatively compact open set
$\Om\Subset X$ and its preimage $\tO=\pi^{-1}\Om$; $\g$ acts on $\tO$ and $\tO/\g=\Om$.
In general we will decorate by tildes the preimages of objects living on the quotient.
Let $U$ be a fundamental domain of the action of $\g$ on $\tO$. This means that (see e.g. \cite{At}):
a) $\tO$ is covered by the translations of $\overline U$, b) different translations of $U$ have empty 
intersection and c) $\overline U\smallsetminus U$ has zero measure (since $\partial\Om$ is smooth).
Since $\Om$ is relatively compact $U$ has the same property.
Let us define the space of square integrable sections $L^2(\tO,\tf)$ with respect to a
$\g$ invariant metric on $\tm$ (and its volume form) and a $\g$ invariant metric on $\tf$.
Then $L^2(U,\tf)$ is constructed with respect to the same. There is a unitary action of 
$\g$ on $L^2(\tO,\tf)$. In fact it is easy to see that
$L^2(\tO,\tf)\cong L^2\g\otimes L^2(U,\tf)\cong L^2\g\otimes L^2(\Om,F)$.

\noindent
We have a  unitary action of $\g$ on $L^2\g$ by left translations: $\gamma\longmapsto
l_\gamma$ where $l_\gamma f(x)=f(\gamma^{-1}x)$ for $x\in\g$, $f\in L^2\g$. It induces 
an action on $L^2(\tO,\tf)$ by $\gamma\longmapsto L_\gamma=l_\gamma\otimes\operatorname{Id}$.
Finally we denote by $\De(.\,,.)$ the various spaces of smooth compactly
supported sections.

Let us consider a formally self-adjoint, strongly elliptic, positive differential operator
$P$ on $M$ acting on sections of $F$. Denote by $\tp$ the $\g$--invariant differential operator 
which is its pull-back to $\tm$. From $\tp$ we construct the following operators: the 
Friedrichs extension in $L^2(\tO,\tf)$ of $\tp$ with domain $\De(\tO,\tf)$ and the Friedrichs extension in
$L^2(U,\tf)$ of $\tp$ with domain $\De(U,\tf)$. From now on we denote these extensions 
by $\tp$ and $\p0$. They are closed self-adjoint positive operators.
It is known that $\tp$
is also $\g$ invariant i.e. it commutes with all $L_\gamma$. This amounts of saying that 
$E_{\lambda}$ commutes with $L_\gamma$, $\gamma\in\g$, where 
$(E_{\lambda})_{\lambda}$ is the spectral family of $\tp$. 
On the other hand
the Rellich lemma tells that $\p0$ has compact resolvent and hence discrete spectrum.
We will take the task of comparing the distribution of the two spectra.  
Namely since $E_\lam$ is $\g$ invariant its image $R(E_\lam)$
is a $\g$ invariant closed subspace of the free Hilbert $\g$--module $L^2\g\otimes L^2(U,\tf)
\cong L^2(\tO,\tf)$. In general for any Hilbert space $\Cal H$ we call the Hilbert space
$L^2\g\otimes\Cal H$ a free Hilbert $\g$--module. The action of $\g$ is defined as above by
$\gamma\longmapsto L_\gamma=l_\gamma\otimes\operatorname{Id}$.
For $\g$ invariant closed spaces 
(called $\g$ modules) 
one can associate a positive, 
possibly infinite real number, called von Neumann or $\g$--dimension, denoted 
$\dim_{\scriptscriptstyle\Gamma}$. For notions involving the $\g$--dimension
and linear algebra for $\g$--modules we refer the reader to \cite{At}, \cite{Sh} and 
\cite{Ko} (in the latter proofs from scratch are given).
We give here the barest discussion of this score. Let us denote by ${\Cal A}_{\sG}$ the 
von Neumann algebra which 
consists of all bounded linear operators in $L^2\g\otimes\Cal H$ which commute to the action of $\g$.
To describe ${\Cal A}_{\sG}$ let us consider the von Neumann ${\Cal R}_{\sG}$ algebra of all 
bounded operators on
$L^2\g$ which commute with all $L_\gamma$. It is generated by all right translations. If we consider 
the orthonormal basis $(\delta_\gamma)_\gamma$ in $L^2\g$ where $\delta_\gamma$ is the Dirac delta 
function at $\gamma$, then the matrix of any operator $A\in{\Cal R}_{\sG}$ has the property that all
its diagonal elements are equal. Therefore we define a natural trace on ${\Cal R}_{\sG}$ as 
the diagonal element, that is,
$\operatorname{tr}_{\sG} A= (A\delta_e,\delta_e)$
where $e$ is the neutral element. Now ${\Cal A}_{\sG}$ is the tensor product of 
${\Cal R}_{\sG}$ and the 
algebra ${\Cal B}(\Cal H)$ of all bounded operators on $\Cal H$. 
If $\operatorname{Tr}$ is the usual trace on ${\Cal B}(\Cal H)$ then we have a trace on
${\Cal A}_{\sG}$ by $\operatorname{Tr}_{\sG} = \operatorname{tr}_{\sG} \otimes\operatorname{Tr}$.
For any $\g$ invariant space $L\subset L^2\g\otimes\Cal H$ i.e. for any $\g$--module, the projection 
$P_L\in{\Cal A}_{\sG}$ and we define $\dg L=\operatorname{Tr}_{\sG}P_L$. 
Let us just remark for
later use that if $L\subset L^2(\tO,\tf)$ is a $\g$--module and
$f_i$ is an orthonormal basis of $L$ then 
$$\dim_{\sG}L=\sum_i\int_{U}|f_i|^2\,.\eqno(1.1)$$

\noindent
We denote in the sequel
$\nlp=\dg R(E_\lam)$. Similary we consider the spectral distribution function 
(counting function) $\np0=\dim R(E^0_\lam)$
where $E^0_\lam$ is the spectral family of $P_0$; it equals the number of eigenvalues 
$\leqslant\lam$. We want to compare $\nlp$ and $\np0$. For this purpose we use essentially 
the analysis of Shubin \cite{Sh}. However there exist a difference in our method, namely we work 
at the beginning with model operator $P_0$ the operator $\tp$ itself with Dirichlet boundary 
conditions on $U$ whereas Shubin considers a direct sum of tangent operators to $\tp$. 
So we do not have to truncate
from the outset the eigenfunctions of the model $P_0$. (See also Remark 1.3 in \cite{Sh} and 
compare e.g. formulas (2.7), (2.8) or (3.6) from \cite{Sh} with our corresponding formulas.)
To begin with we need a variational principle. 

\proclaim{Proposition 1.1(\cite{Sh})}
Let $\tp$ be a $\g$ invariant self-adjoint positive operator on a free $\g$--module
$L^2\g\otimes\Cal H$ where $\Cal H$ is Hilbert space. Then
$$
\eqalignno{
\nlp=\sup\Big\lbrace\dim_{\scriptscriptstyle \Gamma}L\mid
L \;\text{is a}\; \Gamma-\text{module}
\subset&\Dom(\tq),\cr
\tq(f,f)&\leqslant \lambda \|f\|^2,\,\forall f\in L\Big\rbrace 
\quad (1.2)}
$$
where $\tq$ is the quadratic form of $\tp$. 
\endproclaim

\noindent 
Recall that $\tq$ is the closed symmetric
quadratic given by $\Dom(\tq)=\Dom(\tp^{1/2})$, $\tq(u)
=(\tp^{1/2}u,\tp^{1/2}u)$. From the variational principle we deduce the following.

\proclaim{Proposition 1.2 (Estimate from below)} 
The counting functions for $\tp$ and $P_0$ satisfy the inequality

$$\nlp\geqslant\np0\,,\quad \lam\in\Bbb R\eqno(1.3)$$

\endproclaim

\demo{Proof}                                     
Let us denote by $\lam_0\leqslant\lam_1\leqslant\dotsc$ the spectrum of $P_0$. 
Let $\{e_i\}_i$ be an orthonormal basis of
$L^2(U,\widetilde F)$ which consists of eigenfunctions of $P_0$ 
corresponding to the eigenvalues $\{\lambda _i\}_i$; if we let 
$\widetilde e_i=0$ on $\tO\smallsetminus\overline U$ and $\widetilde  
e_i=e_i$ on $U$, $\widetilde e_i\in\Dom(\tq)$,
$\{L_{\gamma}\widetilde e_i\}_{i,\gamma}$
is an orthonormal basis of $L^2(\tO,\widetilde F)$ and $\widetilde
e_{i,\gamma}=L_{\gamma}\widetilde  e_i\in 
\Dom(\tq)$. We have $\tq(\widetilde e_{i,\gamma},\widetilde
e_{i^{'},\gamma^{'}})=\delta_{i,i^{'}}\delta_{\gamma,\gamma^{'}}\lambda_i$.
Let $\Phi_{\lambda}^{0}$
be the subspace spanned by $\{e_i:{\lambda_i\leqslant \lambda}\}$ in $L^2(U,\tf)$ and
$\Phi_{\lambda}$ the closed subspace spanned by $\{\widetilde e_{i,\gamma}: 
{\lambda_i\leqslant \lambda}\}$ in $L^2(\tO,\tf)$. Then by (1.1)
$$
\dim_{\scriptscriptstyle\Gamma}\Phi_{\lambda}=\sum
\limits_{\lambda_i\leqslant\lambda\,\gamma\in\g}
\int_U |\widetilde e_{i,\gamma}|^2
=\sum\limits_{\lambda_i\leqslant\lambda}
\|e_i\|^2_U=\dim\Phi_{\lambda}^0=N(\lambda,P_0)
$$
since $\widetilde e_{i,\gamma}|_U$ vanishes unless $\gamma$ is the identity,
and then it equals $e_i$\,.  
If $f$ is a linear
combination of $\widetilde e_{i,\gamma},\lambda_i \leqslant \lambda$, then
$\tq(f,f)
\leqslant\| f\|^2$ and, as $\Dom(\tq)$ is complete in the graph norm, we obtain that 
$\Phi_{\lambda}\subset\Dom(\tq)$ and $\tq(f,f)\leqslant\lambda\|f\|^2$, 
$f\in \Phi_{\lambda}$. From the variational principle  it follows that 
$N_{\scriptscriptstyle\Gamma}(\lambda ,\tp)\geqslant N(\lambda ,P_0)$.
\hfill{\,}\qed
\enddemo

The next step is an estimate from above of $\nlp$. 
Before let us say something about $\g$--morphisms. If $L_1$, $L_2$
are two $\g$--modules then an bounded liniar operator $T:L_1\longrightarrow L_2$
is called a $\g$--morphism if it commutes with the action of $\g$. As for the usual 
dimension the following statements are true (see \cite{Ko}). If $T$ is injective then
$\dg L_1\leqslant\dg L_2$ and if $T$ has dense image then $\dg L_1\geqslant\dg L_2$.
We denote by $\rank_{\scriptscriptstyle\Gamma}T=\dg\overline{R(T)}$.
For the following we refer to \cite{Sh}, Lemma 3.7.
 
\proclaim{Lemma 1.3}Let us consider the same setting as in the variational principle.
Assume there is $T: L^2(\tO,\widetilde F)\rightarrow L^2(\tO,\widetilde F)$
a $\Gamma$--morphism   
such that
$((\tp+T)f,f)\geqslant\mu \|f\|^2$, $f\in\Dom(\tp)$ and
$\rank_{\scriptscriptstyle\Gamma}T\leqslant p$. Then
$$
N_{\scriptscriptstyle\Gamma}(\mu -\varepsilon ,\tp)\leqslant p,\;\forall 
\varepsilon >0. \eqno(1.4)
$$
\endproclaim

In order to get an estimate from above we have to enlarge a little bit the fundamental 
domain $U$ and compare the counting function of $\tp$ to the counting function of the 
Friedrichs extension of $\tp$ restricted to compactly supported forms in the enlarged domain.
For $h>0$, the enlarged domain is $U_h=\{x\in \tO\mid d(x,U)<h\}$ where $d$ is the 
distance on $\tm$ associated to the Riemann metric on $\tm$. Then we take the tranlations
$U_{h,\gamma}:=\gamma U_h$. Next we construct a partition of unity. 
Let $\varphi^{(h)}\in C^{\infty}(\tO)$, $\varphi^{(h)}\ge 0$, 
$\varphi^{(h)}=1$ on $\bar U$ and $\supp\varphi^{(h)}\subset U_h$,
$\varphi^{(h)}_{\gamma}=\varphi^{(h)}\circ \gamma^{-1}$. 
We define the function $J^{(h)}_{\gamma}\in C^{\infty}(\tO)$ by 
$J^{(h)}_{\gamma}=\varphi_{\gamma}^{(h)}
 \big(\sum_{\gamma}(\varphi^{(h)}_{\gamma})^2\big)^{-\frac{1}{2}}$ 
so that $\sum_{\gamma\in\Gamma}(J^{(h)}_{\gamma})^2=1$.
If $\tp$ is of {\it order $2$}, which will be assumed throughout the section,
then by \cite{Sh,Lemma 3.1} (Shubin's IMS localization formula, see \cite{CFKS})
we know how to recover 
the operator $\tp$ from its localisations $J_{\gamma}^{(h)}\tp J^{(h)}_{\gamma}$:
$$
\tp= \sum\limits_{\gamma \in \Gamma}J_{\gamma}^{(h)}\tp J^{(h)}_{\gamma}-
\sum\limits_{\gamma \in \Gamma}\sigma_0(\tp)(dJ_{\gamma}^{(h)})\eqno(1.5)
$$
where $\sigma_0$ is the principal symbol of $\tp$. 
In (1.5) $J_{\gamma}^{(h)}$ are thought as multiplication operators on 
$L^2(\tO,\tf)$ -- for which $\Dom(\tp)$ is invariant -- while 
$\sum_{\gamma \in \Gamma}\sigma_0(\tp)(dJ_{\gamma}^{(h)})$  
is the multiplication by a bounded function. Since the derivative of 
$J_{\gamma}^{(h)}$ is $O(h^{-1})$ and the order of $\tp$ is $2$ we see that 
the latter function is bounded by $C\,h^{-2}$ for some constant $C>0$ 
(here we use that the symbol is periodic and that $\varphi_{\gamma}^{(h)}$
are the translates of $\varphi^{(h)}$). Therefore the operatorial norm of 
the multiplication satisfies the same estimate and we deduce from (1.5) that
$$
 \tp\geqslant \sum\limits_{\gamma\in\Gamma}J^{(h)}_{\gamma}\tp 
J^{(h)}_{\gamma}-\frac{C}{h^2}\operatorname{Id}\eqno(1.6)
$$
We consider now the operator $\tp$ with domain ${\Cal D}(U_h,\tf)$
and take its Friedrichs extension denoted $P_0^{(h)}$. 
We will compare $\nlp$ with the counting function of $P_0^{(h)}$.
Let us fix $\lam$. Denote by $(E_\lam^{(h)})$ the spectral 
family of $P_0^{(h)}$ and fix a positive constant $M=M^{(\lam)}$ such that
$M\geqslant \lam - \inf\operatorname{spectrum}P_0^{(h)}$
to the effect that
$P_0^{(h)} + M\,E_\lam^{(h)}\geqslant\lam\operatorname{Id}$\,.
We define now a localisation of $E^{(h)}_{\lam}$ by taking the bounded operators
$G^{(h)}_{\gamma}$ on $L^2(\tO,\tf)$ given by 
$G^{(h)}_{\gamma}=J^{(h)}_{\gamma}\,L_{\gamma}\,M\,E^{(h)}_{\lam}
\,L_{\gamma}^{-1}\,J_{\gamma}^{(h)}$ and then summing over $\g$,
$G^{(h)}=\sum_{\gamma\in\Gamma}G^{(h)}_{\gamma}$.
We have
$$
\eqalignno{
\tp+G^{(h)}&\geqslant\sum\limits_{\gamma\in\Gamma}\left(
J^{(h)}_{\gamma}\tp J^{(h)}_{\gamma}
+J^{(h)}_{\gamma}L_{\gamma}M\,E^{(h)}_{\lambda}L_{\gamma}^{-1}
J_{\gamma}^{(h)}\right)-\frac{C}{h^2}\operatorname{Id}\cr
&=\sum\limits_{\gamma\in\Gamma}J_{\gamma}^{(h)}L_{\gamma}(H^{(h)}_0+ 
M\,
E^{(h)}_{\lambda})L_{\gamma}^{-1}J_{\gamma}^{(h)}-\frac{C}{h^2}
\operatorname{Id}\qquad\qquad\qquad\qquad\quad(1.7)\cr%
&\geqslant\sum\limits_{\gamma\in\Gamma}J_{\gamma}^{(h)}L_{\gamma}\lambda 
L_{\gamma}^{-1}
J_{\gamma}^{(h)}-\frac{C}{h^2}\operatorname{Id}
=\left(\lambda -\frac{C}{h^2}\right)\operatorname{Id}.
}
$$
It is clear that $G^{(h)}$ will play the role of $T$ in Lemma 2.3. 
We must check one more hypothesis.

\proclaim{Claim 1.4}
$$\rank_{\scriptscriptstyle\Gamma}G^{(h)}\le N(\lambda,P^{(h)}_0)\eqno(1.8)$$ 
\endproclaim

\demo{Proof}  We start with 
the bounded operator
$\bar G^{(h)}$ on $L^2(U_s,\tf)$, given by
$\bar G^{(h)}=J_{e}^{(h)}M\,E_{\lambda}^{(h)}J_{e}^{(h)}$\,.
It is a finite rank operator, $\rank\bar G^{(h)}\leqslant\rank E_\lam^{(h)}
=N(\lam,P_0^{(h)})$. Next we consider the free $\g$--module $L^2\g\otimes
L^2(U_h,\tf)$ and the bounded $\g$--invariant operator $\Id\otimes\bar G^{(h)}$.
Then $R(\Id\otimes\bar G^{(h)})=L^2\g\otimes R(\bar G^{(h)})$ so that
$\rank_{\sG}\Id\otimes\bar G^{(h)}=\rank\bar G^{(h)}$. We identify now the
space $L^2\g\otimes L^2(U_h,\tf)$ with $\bigoplus_{\gamma\in\Gamma}
L^2(U_{h,\gamma},\tf)$ by the unitary transform 
$K:\,\sum_{\gamma}\delta_{\gamma}\otimes w_\gamma\longmapsto 
\left(L_\gamma w_\gamma\right)_\gamma$\,.
Thus $\bigoplus_{\gamma\in\Gamma}L^2(U_{h,\gamma},\tf)$ is naturally
a free $\g$--module for which $K$ is $\g$ invariant.
We transport $\Id\otimes\bar G^{(h)}$ on $\bigoplus_{\gamma\in\Gamma}L^2(U_{h,\gamma},\tf)$ 
by $K$ and we think it as acting on this latter space.
We construct then a restriction operator
$
V: \bigoplus_{\gamma\in\Gamma} L^2(U_{h,\gamma},\tf)\longrightarrow L^2(\tO,\tf)$\,,
$V\left((w_{\gamma})_{\gamma}\right) =\sum_{\gamma\in\Gamma}w_{\gamma}$\, 
which is a surjective $\g$--morphism. We have also the $\g$--morphism
$I$ from $L^2(\tO,\tf)$ to $\bigoplus_{\gamma\in\Gamma}L^2(U_{h,\gamma},\tf)$,
$I(u)=(u\upharpoonright_{U_{h,\gamma}})_{\gamma}$\, which is obviously bounded. 
With our identifications we have
$G^{(h)}=V\,(\Id\otimes\bar G^{(h)})\,I$\,.
As in the case of usual dimension 
$\rank_{\sG} V\,(\Id\otimes\bar G^{(h)})\,I\leqslant\rank_{\sG}(\Id\otimes\bar G^{(h)})$
(see \cite{Sh},
Lemma 3.6).Therfore we conclude 
$\rank_{\sG} G^{(h)}\leqslant
\rank_{\sG}(\Id\otimes\bar G^{(h)})=
\rank \bar G^{(h)}\leqslant N(\lam,P_0^{(h)})$\,.
\hfill{\,}\qed
\enddemo

\proclaim{Proposition 1.5 (Estimate from above)}  
There is a constant $C\ge 0$ such that 
$$
N_{\sG}(\lambda ,\tp)\leqslant N\left(\lambda 
 + \frac{C}{h^2},P_0^{(h)}\right)
\quad\lambda\in{\Bbb R},\quad h>0
\eqno(1.9)$$
\endproclaim 

\demo{Proof} The hypothesis of Lemma 2.3 are fulfilled for $T=G^{(h)}$, $\mu=\lam-C\,h^{-2}$
and $p=N(\lam,P_0^{(h)})$ as (1.7) and (1.8) show. Thus
$N_{\scriptscriptstyle\Gamma}
\left(\lambda-\frac{C}{h^2}-\varepsilon ,\tp\right)
\leqslant N(\lambda , P_0^{(h)})$, if $\varepsilon >0$.
Replacing $\lambda$ with $\lambda+C\,h^{-2}+\varepsilon$,
we obtain
$\nlp\leqslant N\left(\lambda +\frac{C}{h^2}+\varepsilon ,P_0^{(h)} \right)$.
When $\varepsilon \longrightarrow 0$ the estimate (1.9) follows since
the spectrum distribution function is right continuous by definition.
\hfill{\,}\qed
\enddemo

The estimates from below and above for $\nlp$ enable us to study as a by--product
the behaviour for $\lam\longrightarrow\infty$ to obtain the Weyl asymptotics for 
periodic operators (Shubin, see \cite{RSS} and the references therein).
\proclaim{Corollary 1.6}
If $\tp$ is a periodic, positive, second order elliptic operator as above then
$$\eqalign{\lim_{\lam\rightarrow\infty}\lam^{-n/2}\nlp &=
\lim_{\lam\rightarrow\infty}\lam^{-n/2}N(\lam, P_0)\cr
&=(2\pi)^{-n}\int_U \int_{T^*_{x}\tm}N(1,\sigma_0(\tp)(x,\xi))d\xi\,dx}
$$
where $\sigma_0(\tp)(x,\xi)\in\operatorname{Herm}(\tf,\tf)$ is the 
principal symbol of $\tp$ and  $N(1,\sigma_0(\tp)(x,\xi))$ is the counting function 
for the eigenvalues of this hermitian matrix.
\endproclaim
\demo{Proof}
First let us remark that the last equality is the classical Weyl type formula as 
established by Carleman, G\aa rding and others, see \cite{RSS}, p.72.
It is obvious that $\liminf\lam^{-n/2}N(\lam, P_0)\leqslant\liminf\lam^{-n/2}\nlp$
by the estimate from below. On the other hand the estimate from above gives
$$\eqalign{
\limsup\lam^{-n/2}\nlp\leqslant\limsup\left(1+\frac{C}{\lam\,h^2}\right)^{n/2}
\left(\lam +\frac{C}{h^2}\right)^{-n/2} N\left(\lambda 
 + \frac{C}{h^2},P_0^{(h)}\right)\cr
\leqslant\limsup\mu^{-n/2}N(\mu,P_0^{(h)})=
(2\pi)^{-n}\int_{U_h} \int_{T^*_{x}\tm}N(1,\sigma_0(\tp)(x,\xi))d\xi\,dx
}$$
for a fixed small $h$. We make $h\longrightarrow 0$ and obtain the desired
formula.
\hfill{\,}\qed
\enddemo

We are going to apply the above results to the semi-classical asymptotics
as $k\longrightarrow\infty$ of the spectral distribution function of the laplacian
$k^{-1}\lat$ on $\tm$.
Let $G$ be a hermitian holomorphic bundle on $M$ and
$\widetilde G=p^\ast G$ its pull-back.
We define ${\Cal D}^{(0,q)}(.\,,.)$ to be the space of smooth compactly supported $(0,q)$ forms.
Let $\bar\partial:{\Cal D}^{0,q}(\tm,\widetilde G)\longrightarrow
{\Cal D}^{0,q+1}(\tm,\widetilde G)$ 
be the Cauchy--Riemann operator and $\vartheta:{\Cal D}^{0,q+1}(\tm,\widetilde G)
\longrightarrow {\Cal D}^{0,q}(\tm,\widetilde G)$ 
the formal adjoint of $\bar\partial$ with respect to the given hermitian metrics on $\tm$, 
$\tg$. Then $\lat=\bar\partial\vartheta +\vartheta\bar 
\partial$ is a formally self-adjoint, strongly elliptic, positive and $\g$--invariant 
differential operator.

We take $\te$ and $\tg$ two $\g$ invariant holomorphic bundles. 
Let us form the 
Laplace--Beltrami operator $\latk$ on $(0,q)$ forms with values in $\te^k\otimes\tg$.
Thus we will consider the $\g$ invariant hermitian bundle $\tf=\Lambda^{(0,q)}T^{\ast}
\tm\otimes\te^k\otimes\tg$ and apply the previous results for $\tp=k^{-1}\latd$ where the 
index $\tO$ emphasises that the Friedrichs extension gives the operator of the Dirichlet problem
on $\tO$.
Now we have to make a good choice of the parameter $h$. We take $h=k^{-\frac{1}{4}}$
so that the derivative of the cutting off function $J_\gamma^{(h)}$ is just 
$O(k^{\frac{1}{4}})$. Then $\sigma_0(k^{-1}\latk)(dJ_{\gamma}^{(h)})=k^{-1}|\bar\partial 
J_{\gamma}^{(h)}|^2=O(k^{-\frac{1}{2}})$. Therefore formula (1.6) becomes
$
 \frac{1}{k}\latd\geqslant \sum_{\gamma\in\Gamma}J^{(h)}_{\gamma}\frac{1}{k}\latd 
J^{(h)}_{\gamma}-\frac{C}{\sqrt{k}}\operatorname{Id}\,.
$  
We have thus proved the following semi--classical estimate for laplacian.

\proclaim{Proposition 1.7}
There exists a constant $C>0$ such that for $\lambda\in{\Bbb R}$ and $ k>0$ we have
$$
N\left(\lam\,,\frac{1}{k}\latk\upharpoonright _U\right)\leqslant
N_{\sG}\left(\lambda\,,\frac{1}{k}\latd\right)\leqslant 
N\left(\lambda  + \frac{C}{\sqrt{k}}\,,\frac{1}{k}\latk
\upharpoonright _{ U_{k^{-1/4}}}\right)
\eqno(1.10)
$$
\endproclaim

Demailly has determined the distribution of spectrum for the Dirichlet problem for
$\latk$ in \cite{De1}, Theorem 3.14. For this purpose he introduces (\cite{De1},(1.5))
the function $\nu_E :\tm\times\Bbb R\longrightarrow\Bbb R$ depending on the curvature of $\te$
and then considers the function $\bar\nu_E(x,\lam)=\lim_{\ve\searrow 0}\nu_E(x,\lam+\ve)$.
The function $\bar\nu_E(x,\lam)$ is right continuous in $\lam$ and bounded above 
on compacts of $\tm$.
Denote by $\alpha_1(x),\dotsc,\alpha_n(x)$ the eigenvalues of of the 
curvature form $i{\bold c}(\te)(x)$ with respect to the metric on $\tm$. 
We also denote for a multiindex $J\subset\{1,...,n\}$,
$\alpha_J=\sum_{j\in J}\alpha_j$ and $C(J)=\{1,...,n\}\smallsetminus J$\,. 
For $V\Subset M$ we introduce
$$I^q(V,\mu)=\sum\limits_{\mid J\mid =q}\int_{V}\bar\nu_E(2\mu+\alpha_{C(J)}-\alpha
_{J})\,d\sigma$$

\proclaim{Proposition 1.8 (Demailly)}
Assume that $\partial V$ has measure zero and that the laplacian acts on $(0,q)$ forms.
Then
$\limsup_{k}\,k^{-n}N(\lambda\,,\frac{1}{k}\Delta^{\prime\prime}_{k}
\upharpoonright _V)\leqslant I^q(V,\lam)
$ 
Moreover there exists an at most countable set
$\Cal N\subset\Bbb R$ such that for $\lam\in\Bbb R\smallsetminus\Cal N$ the limit of the left--hand
side expression exists and we have equality.
\endproclaim

We return now to the case of a covering manifold and apply Demailly's formula 
in (1.10).
Let us fix $\ve >0$. For sufficiently large $k$ we have $U_{k^{-\frac{1}{4}}}
\subset U_\ve$
so the fact that the counting function is increasing and the variational principle yield
$
N(\lambda+\frac{C}{\sqrt{k}}\;,\frac{1}{k}\lat
\upharpoonright _ {U_{k^{-1/4}}})\leqslant
N(\lambda+\ve,\frac{1}{k}\lat
\upharpoonright _ {U_{k^{-1/4}}})\leqslant
N(\lambda+\varepsilon 
,\frac{1}{k}\lat\upharpoonright _{U_{\varepsilon}})
$.
Hence by (1.10) and Proposition 1.8 ($\partial U_\ve$ is negligible for small $\ve$),
$$\limsup_{k} k^{-n}N_{\sG}(\lambda\;,\frac{1}{k}\latd)\leqslant I^q(U_\ve,\lam+\ve).$$
The use of dominated convergence to make 
$\ve\longrightarrow 0$ in the last integral yield
the asymptotic formula for the laplacian on a covering manifold.
\proclaim{Theorem 1.9}
The spectral distribution function of $\frac{1}{k}\latd$
on $L^2_{0,q}(\widetilde\Om,\te^k\otimes\tg)$ with Dirichlet boundary conditions satisfies
$$\limsup_{k}\,k^{-n}N_{\sG}\left(\lambda\,,\frac{1}{k}\latd\right) \leqslant
I^q(U,\lam).\eqno(1.11)
$$
Moreover, there exists an at most countable set 
$\Cal N\subset\Bbb R$ such that for $\lam\in\Bbb R\smallsetminus\Cal N$ the limit exits and
we have equality in \rom{(1.11)}.
\endproclaim

\subhead
\S 2 Geometric situations
\endsubhead

In this section we apply the results from the previous section to the study
of the $L^2$ cohomology of coverings of complex manifolds satisfying certain curvature conditions.
If $M$ is a complete K\"ahler manifold and $E$ a positive line bundle on $M$ the $L^2$ estimates of
Andreotti--Vesentini--H\"ormander allow to find a lot of sections of $\te$ on a covering $\tm$ (see e.g. 
\cite{NR}). We prove here the following.
\proclaim{Theorem 2.1}
Let $(M,\omega)$ be an $n$--dimensional complete hermitian manifold such that the torsion
of $\omega$ is bounded 
and let $(E,h)$ be a holomorphic hermitian line bundle.
Let $K\Subset M$ and a constant $C_0>0$ such that $\imath{\boldkey c}(E,h)\geqslant C_0\omega$
on $M\smallsetminus K$.
Let $p:\tm\longrightarrow M$ be a Galois covering with group
$\g$ and $\te=p^\ast E$ and let $\Om$ be any open set with smooth boundary and $K\Subset\Om\Subset M$.
Then 
$$\dg H^{n,0}_{(2)}(\tm,\te^k)\geqslant \frac{k^n}{n!}\int_{\Om(\leqslant 1,h)}
\left(\frac{\imath}{2\pi}{\boldkey c}(E, h)\right)^n +o(k^n)
\,,\quad k>>0\,, $$
where $H^{n,0}_{(2)}(\tm,\te^k)$ is the space of $(n,0)$--forms with values in
$\te^k$ which are $L^2$ with respect to any metric on $\tm$ and the pullback of $h$ 
and $\Om(\leqslant 1,h)$ is the subset of $\Om$ where $\imath{\boldkey c}(E,h)$ is
non--degenerate and has at most one negative eigenvalue.
\endproclaim

\demo{Proof}
We endow $\tm$ with the metric $\widetilde\omega=p^{\ast}\omega$ and $\te$ with 
$\widetilde h=p^{\ast}h$. All the norms, Laplace--Beltrami operators,
spaces of harmonic forms and $L^2$--cohomology groups are with respect to $\widetilde\omega$
and $\widetilde h$. In particular the operators $\db$ and Lapalce--Beltrami are $\g$--invariant. 
It is standard to see that $\widetilde\omega$ is also complete. 
To justify this 
let us first take a compact set $K\Subset X$ and consider $\widetilde K=p^{-1}K$.
The metric $\widetilde\omega_\ve$ is complete on $\widetilde K$ in the following sense.
There exist functions $\varphi_\ve\in{\Cal C}^\infty (\widetilde K)$ with values in $[0,1]$
such that $\operatorname{supp}\varphi_\ve$ is compact in $\widetilde K$, the sets 
$\{z\in\widetilde K\,:\,\varphi_\ve(z)=1\}$ form an exhaustion of $\widetilde K$
and $\sup |d\varphi_\ve|=O(\ve)$ as $\ve\longrightarrow 0$. This is seen as usual by 
observing that the balls are relatively compact in $\widetilde K$
and then taking cut--off functions. Since $M$ is complete there exist an exhaustion
$K_\nu$ with compacts and functions $\psi_\nu\in{\Cal C}^\infty(M)$ with values in $[0,1]$ and
$\operatorname{supp}\psi_\nu\Subset K_{\nu +1}$
such that $K_\nu\subset\{z\in M\,:\,\psi_\nu(z)=1\}$ and $\sup |d\psi_\nu|\leqslant
2^{-\nu}$. Let us choose now a point $z_0\in\widetilde K_0$ and fix fundamental
domains $U_\nu$ for the action of $\g$ on $\widetilde K_\nu$ such that $z_0\in U_\nu$. 
We also choose an exhaustion by finite sets $I_0\subset I_1\subset\dotsm\subset 
I_\nu\subset\dotsm\subset\g$ of $\g$. Indeed, since $\tm$ is paracompact $\g$ is countable. 
For each $\nu$ let us take $\varphi_\nu\in{\Cal C}^\infty (\widetilde K_{\nu+1})$ such that 
$\varphi_\nu=1$ on $\cup\{\gamma U_{\nu+1}\,:\,\gamma\in I_{\nu+1}\}$ and 
$\sup |d\varphi_\nu|\leqslant 2^{-\nu}$. We consider also the function $\widetilde\psi_\nu
=\psi_\nu\circ p$. Then the functions  $\widetilde\psi_\nu\varphi_\nu$ have compact support
in $\tm$, the sets where they equal $1$ exhaust $\tm$ and their derivative is $O(2^{-\nu-1})$,
which proves that $\tm$ is complete.
We remark here that ${\Cal U}=\cup_{\nu}U_{\nu}$ is a fundamental domain for the action of $\g$ on
$\tm$ and that if $\widetilde G$ is a $\g$--invariant bundle on $\tm$ then $L^2(\tm,
\widetilde G)$ is a free $\g$--module.

We take $\Om$ as in the hypothesis
and let $U$ be a fundamental domain of $\tO$ as in \S 1.
Since $p$ is locally biholomorphic we see that 
$\imath{\bold c}(\tE,\widetilde h)\geqslant C_0\widetilde\omega$ 
on $\tm\smallsetminus\widetilde K$.
Let $u$ be a smooth $(n,1)$ form on $\tm$
with values in $\te^k$ and compactly supported outside
$\widetilde K$. 
We apply now the Bochner--Kodaira--Nakano formula for $u$:
$$3\left(\latk u,u\right)\geqslant 
2\left(\left[\imath{\boldkey c}(\te^k),\widetilde\Lambda\right]u,u\right)
-\left(\|\tau\,u\|^2+\|\bar\tau\,u\|^2
+ \|\tau^{\ast}\,u\|^2+\|\bar\tau^{\ast}\,u\|^2\right)\,,
$$
where $\Lambda$ is the operator of taking the interior product with
$\widetilde\omega$ and the $\tau$'s are the torsion operators of the
metric $\widetilde\omega$. More precisely $\tau=[\Lambda,\partial\widetilde\omega]$.
Therefore there exists a constant $C_1>0$ (depending just on the metric $\omega$)
such that
$$3\left(\latk u,u\right)\geqslant
2C_0\,k\,\|u\|^2 - C_1\,\|u\|^2\,,$$
and hence
$$\left(\latk u,u \right)\geqslant\frac{C_{\scriptscriptstyle 0}\,k}{2}\;\|u\|^2\,,\quad
k\geqslant \frac{\,C_1}{2C_{\scriptscriptstyle 0}}\,.\eqno(2.1)$$
Indeed, by hypothesis the torsion operators are pointwise bounded. 
Moreover
$([\imath{\boldkey c}(\te^k,\widetilde h^k),\Lambda]u,u)\geqslant
k\,\alpha_1\,|u|^2$ where $\alpha_1\leqslant\dotsm\leqslant\alpha_n$
are the eigenvalues of $\imath{\boldkey c}(\te,\widetilde h)$ with respect to
$\widetilde\omega$.

Let $\rho\in{\Cal C}^{\infty}(M)$ such that $\rho=0$
on $L$ and $\rho=1$ on $M\smallsetminus\Om$, where $L$ is a neighbourhood of $K$ in $\Om$. 
We put $\widetilde\rho=\rho\circ p$. Let $u\in\De^{n,1} (\tm ,\te^k)$,
so that $\widetilde\rho\,u$ has support outside $\widetilde K$.
We use now the elementary estimate:
$$\left(\latk (\widetilde\rho\,u) ,\widetilde\rho\,u\right)\leqslant
\frac{3}{2}\left(\latk u,u\right)+6\sup |d\widetilde\rho\,|^2
\|u\|^2\,.\eqno(2.2)
$$
Obviously $C_2=6\sup |d\widetilde\rho\,|^2<\infty$. Estimates (2.1) and
(2.2) yield
$$\|u\|^2\leqslant\frac{12}{C_{\scriptscriptstyle 0}\,k}\left(\latk u,u\right) 
+ 4\int_{\tO}\big|(1-\widetilde\rho\,)u\big|^2\,,\quad
k\geqslant \frac{\max\{C_1,16C_2 \}}{2C_{\scriptscriptstyle 0}}\eqno(2.3)$$
for 
any compactly supported $u$. Since the metric $\widetilde\omega$
is complete the density lemma of Andreotti and Vesentini \cite{AV}
shows that $\latk$ is essentially self--adjoint. Thus (2.3) is true for any 
$u$ in the domain of the quadratic form $\widetilde Q_k$ of the self--adjoint
extension of $k^{-1}\latk$. 
From relation (2.3) we infer that the spectral spaces corresponding to the
lower part of the spectrum of $k^{-1}\latk$ on $(n,1)$--forms can be injected into the spectral spaces  
of the $\g$--invariant operator $k^{-1}\latd$ which correspond to the
Dirichlet problem on $\tO$ for $k^{-1}\latk$.
The latter operator was studied in \S 1.
This idea appears in Witten's proof (see Henniart \cite{He}) and in \cite{Bou} in 
the context of $q$--convex manifolds in the sense of Andreotti--Grauert. 
We claim that for $\lam<C_{\scriptscriptstyle 0}/24$, 
$$L^1_k(\lam)\longrightarrow L^1_{k,\tO} (12\lam+C_3 k^{-1})\,,\quad
u\longmapsto E_{12\lam+C_3 k^{-1}}(k^{-1}\latd) (1-\widetilde\rho)u\,,\eqno(2.4)$$ 
is an injective $\g$--morphism, where $L^1_k(\lam)=\operatorname{Range}\big(E_\lam(k^{-1}\latd)\big)$ 
is the spectral space of $k^{-1}\latk$ on $(n,1)$--forms, $L^1_{k,\tO}(\mu)=
\operatorname{Range}E_{\mu}(k^{-1}\latd)$, the spectral spaces of   
$k^{-1}\latd$ and $C_3=8\,C_2$. 
To prove the claim let us remark that the map (2.4) is the restriction of an operator
on $L^2_{0,1}(\tm,\te^k\otimes K_{\tm})$ 
of the same form; this is continuous and $\g$--invariant
being a composition of a multiplication with a bounded $\g$--invariant function and a 
$\g$--invariant projection. To prove the injectivity we 
choose $u\in L^1_k(\lam)$, $\lam<C_{\scriptscriptstyle 0}/24$ to the effect that
$\widetilde Q_k(u)\leqslant\lam\|u\|^2\leqslant (C_{\scriptscriptstyle 0}/24)\|u\|^2$. Plugging this relation in
(2.3) we get
$$\|u\|^2\leqslant 8\int_{\tO}\big|(1-\widetilde\rho\,)u\big|^2
\,,\quad u\in L^1_k(\lam)\,,\quad\lam<C_{\scriptscriptstyle 0}/24\,.\eqno(2.5)$$
Let us denote by $\widetilde Q_{k,\tO}$ the quadratic form of $k^{-1}\latd$.
Then by (2.2) and (2.5),
$
\widetilde Q_{k,\tO}\big((1-\widetilde\rho)u\big) \leqslant
\frac{3}{2}\,\widetilde Q_{k}(u) + \frac{C_2}{k}\|u\|^2\leqslant
\,\left(12\,\lam +\frac{8\,C_2}{k}\right)\int_{\tO}\big|(1-\widetilde\rho\,)u\big|^2
$
which shows that if $E(12\lam+C_3 k^{-1},k^{-1}\latd)\,(1-\widetilde\rho)u=0$
then $(1-\widetilde\rho)u=0$ so that $u=0$ by (2.5). Therefore (2.4) is injective
and hence
$$N^1_{\sG}\left(\lam,\frac{1}{k}
\latk\right)\leqslant 
N^1_{\sG}\left(12\lam +\frac{C_3}{k},\latd\right)
\,,\quad \lam <(C_{\scriptscriptstyle 0}/24)\,,\eqno(2.6)$$
and thus the spectral spaces $L^1_k(\lam)$, $\lam <C_{\scriptscriptstyle 0}/24$, 
are of finite $\g$--dimension.

Now we can apply Theorem 1.9 for $k^{-1}\latd$ on $\tO$ (with $\tg=K_{\tm}$).
By the variational principle we have that $N^0_{\sG}(\lam,\frac{1}{k}\latk)\geqslant
N^0_{\sG}(\lam,\frac{1}{k}\latd)$ and by Theorem 1.9 for $q=0$
$$\liminf_k k^{-n}N^0_{\sG}\left(\lam,\frac{1}{k}\latk\right)
\geqslant I^0(U,\lam)\,,\quad\lam<C_0/24\,,\quad\lam\in{\Bbb R}\smallsetminus{\Cal N}\eqno(2.7)$$
We find now an upper bound. Fix an arbitrary $\delta >0$. For $k>C_3/\delta$ we have
$N^1(\lam,k^{-1}\latk)\leqslant N^1_{\sG}(12\lam +{C_3}\,k^{-1}\latd)\leqslant
N^1_{\sG}(12\lam +\delta,\latd)$ hence by (1.11) 
$\limsup_{k}k^{-n}N^1_{\sG}(\lam,\frac{1}{k}\latk)
\leqslant I^1(U,12\lam+\delta)$. 
We can let $\delta\longrightarrow 0$ so that
$$\limsup_k k^{-n}N^1_{\sG}\left(\lam,\frac{1}{k}\latk\right)
\leqslant I^1(U,12\lam)\,,\quad\lam<C_0/24\,.\eqno(2.8)$$

We consider the group
$H^{n,0}_{(2)}(\tm,\widetilde E^k)=
\lbrace u\in L^2_{n,0}(\tm,\te^k, \widetilde\omega,{\widetilde h})
\,:\,\db u=0\rbrace$ which is a $\g$--module
and we find a lower bound for its $\g$--dimension.
We know that the $L^2$ norm doesn't actually depend on the metric on $\tm$.
We consider also the operator $\latk$ defined on 
$L^2_{n,0}(\tm,\te^k)$ 
and denote by $L^0_k(\lam)$ its spectral spaces.
Since $\latk$ commutes with $\db$ it follows that the spectral projections of 
$\latk$ commute with $\db$ too, showing thus $\db L^0_k(\lam)\subset L^{1}_k(\lam)$
and therefore we have the $\g$--morphism
$L^0_k(\lam)@>{\;\db_\lam\;}>>L^1_k(\lam)$
where $\db_\lam$ denotes the restriction of $\db$ (by the definition of $L^0_k(\lam)$,
$\db_\lam$ is bounded by $k\lam$). 
Since for any $\g$--morphism
$A$ we have $\dg\overline{R(A)}=\dg\ker (A)^{\perp}$ we see that 
$\dg\ker\db_\lam +\dg\overline{R(\db_\lam)}=\dg L^0_k(\lam)$. Moreover $\dg\overline{R(\db_{\lam})}
\leqslant \dg L^1_k(\lam)$ and they are finite. Therefore by (2.7) and (2.8),
$
\dg H^{n,0}_{(2)}(M,\widetilde E^k)\geqslant\dg\ker\db_\lam
\geqslant k^n\Big[ I^0(U,2\lam) -I^1(U,12\lam)\Big]
$
for $\lam<C_0/24$ and $\lam\in\Bbb R\smallsetminus\Cal N$. 
We can now let $\lam$ go to zero through these values. The limits $I^0(U,0)$ and $I^1(U,0)$
are calculated in \cite{De1} and if we identify the fundamental domain $U$ with $\Om$
the result is exactly the integral from the conclusion. 
\hfill{\,}\qed
\enddemo
To state the following result let us
remind that by the definition of Andreotti and Grauert \cite{AG} a
manifold is called $1$--\,concave if there exists a smooth function $\varphi:X
\longrightarrow (a,b\,]$ where $a\in\{-\infty\}\cup\Bbb R$, $b\in 
\Bbb R$, such that $X_c:=\{\f>c\}\Subset X$ for all $c\in (a,b\,]$ 
and $\f$ is strictly plurisubharmonic outside a compact set.
Let $E$ be a holomorphic line bundle on $X$.
In \cite{Oh}, \cite{Ma} one constructs a function $\chi:(-\infty,0)\longrightarrow\Bbb R$
such that $\int_{-1}^0\chi(t)^{1/2}dt=\infty$, $\chi^{\prime}(t)^2\leqslant4\chi(t)^3$ ,
$\chi(t)\geqslant 4$ and
a hermitian metric $\omega$ which equals $\frac{1}{3}\partial\db\f$ near $bX_c$.
For convenience we denote $\psi=c-\f$.
We define $\omega_0=\omega+\chi(\psi)\partial\f\wedge\db\f$, a complete metric on $X_c$ and a 
hermitian metric $h_0=h\exp(-A\int_{\inf\psi}^{\psi}\chi(t)dt)$ on $E$ over $X_c$.
\proclaim{Theorem 2.2}
Let $X$ be a $1$--concave manifold of dimension $n\geqslant 3$ and let $X_c$ be a sublevel set
such that the exhaustion function $\f$ is strictly plurisubharmonic near $bX_c$.
Let $p:\tX_c\longrightarrow X_c$ be a Galois covering of group $\g$. Assume that $\tX_c$
and $\te$ are endowed with the lifts of the metrics $\omega_0$ and $h_0$. Then
$$\dg H^0_{(2)}(\tX_c,\te^k)\geqslant \frac{k^n}{n!}\int_{\Om(\leqslant 1,h_0)}
\left(\frac{\imath}{2\pi}{\boldkey c}(E, h_0)\right)^n +o(k^n)
\,,\quad k>>0\,. \eqno{2.9}$$
for any sufficiently large open set $\Om\Subset X_c$.
\endproclaim
\demo{Proof} The metrics $\omega_0$ and $h_0$ satisfy the following conditions:
\item{(i)} Denoting by $\gamma_i$ the eigenvalues of $\imath\chi(\psi)\partial\db\psi+
\imath\chi^\prime(\psi)\partial\psi\wedge\db\psi$ with respect to
$\omega_{\scriptscriptstyle 0}$ we have 
$\gamma_1\leqslant\cdots\leqslant\gamma_{n-1}\leqslant -2\chi(\psi)$ and 
$\gamma_n\leqslant\chi(\psi)$ so that $\gamma_n+\cdots+\gamma_2\leqslant(5-2n)\chi(\psi)\leqslant
-\chi(\psi)$ for $n\geqslant 3$ outside a compact set $K:=X_e\Subset X_c$.
\item{(ii)} The torsion operators of the metric $\omega_0$ are pointwise bounded by 
$C_2\chi(\f)^{1/2}$ outside $K$.
\item{(iii)} The eigenvalues of $\imath{\boldkey c}(E,h_0)$ with respect to $\omega_0$
are bounded above on $X_c$ by $C_1>0$.

Let us take the lifts $\widetilde\omega_0$, $h_0$ and $\widetilde\psi=c-\f\circ p$.
It is easy to see that properties (i), (ii) and (iii) are still valid for 
$\widetilde\omega_0$ and $\widetilde h_0$ and $\widetilde\psi$ on $\tX_c\smallsetminus\widetilde K$.
For $u\in\De^{(0,1)}(\tX_c\smallsetminus\widetilde K,\te^k)$ we apply the 
Bochner--Kodaira--Nakano inequality and take into account the formula
$([\imath{\boldkey c}(\te^k,\widetilde h^k),\Lambda]u,u)\geqslant
-k(\alpha_n+\dotsc+\alpha_2)|u|^2$. Then 
$$3\left(\latk u,u\right)\geqslant
\int \left(-knC_1+kA\chi(\psi)-4C_2\chi(\psi)\right)|u|^2\,.$$
For sufficiently $A$ and since $\chi\geqslant 4$ we derive easily an estimate analogous to
(2.1). From this point the proof of Theorem 2.1 applies whith just notational changes.
\hfill{\,}\qed
\enddemo

\subhead
\S 3  Coverings of some strongly pseudoconcave manifolds 
\endsubhead

Let us recall the solution 
of the Grauert-Riemenschneider conjecture (\cite{GR}, p. 277) as given by Siu 
\cite{Si} and Demailly \cite{De1}.
Namely the Siu--Demailly criterion says that if
$X$ be a compact complex manifold and $E$ a line bundle over $X$.
and either 
$E$ is semi-positive and positive at one point (Siu's criterion), or
$$
\int_{X(\leqslant 1)}\big(\imath{\boldkey c}(E)\big)^n>0
\eqno(D)$$
(Demailly's criterion)
then $\dim H^0(X,E^k)\approx k^n$, for large $k$ and $X$ is Moishezon. 
Our aim is to extend this result in two directions.
We allow $X$ to belong to certain classes of strongly pseudoconcave manifolds
and we study (directly) Galois coverings of such manifolds.

For $1$--\,concave and compact manifolds (all which are pseudoconcave in the sense of Andreotti
\cite{An}) the transcendence degree of the meromorphic function field is 
less than or equal to the dimension of $X$. In the latter
case we say that the manifold is Moishezon by extending the
terminology from compact manifolds. 

If, in the Andreotti--Grauert definition, the function $\f$ can be taken such that
$a=\inf\f=-\infty$, we say that $X$ is {\it hyper $1$--\,concave\/}.
Let us note that not all $1$--\,concave manifolds are hyper
$1$--\,concave. Indeed, the complement of $S^1\subset\Bbb C\subset
\Bbb P^1$ in $\Bbb P^1$ is $1$--\,concave but cannot possibly 
be hyper $1$--\,concave since $S^1$ is not a polar set in $\Bbb C$ 
(I have learnt this example from M. Col\c{t}oiu and V. V\^aj\^aitu).

Let us describe some examples.
Let $Y$ be a compact complex manifold, $S$ a complete pluripolar set
(the set where a strictly psh function takes the value $-\infty$).
Then $M=Y\smallsetminus S$ is hyper 1--\,concave. Conversely, if $\dim M\geqslant 3$ 
any hyper 1--\,concave manifold $M$ is biholomorphic to a complement of a pluripolar set
in a compact manifold as a consequence of Rossi's compactification theorem.
Another example of hyper $1$--\,concave manifold is $\rg$ where $M$ is a
compact complex space with isolated singularities.
Suppose that $p$ is an isolated singular point and that the germ
$(X,p)$ is embedded in the germ $({\Bbb C}^N, 0)$ and
$z=(z_1,\dots ,z_N)$ are local coordinates in the ambient space ${\Bbb C}^N$.
The function $\f$ is then obtained by cutting-off functions of the
type $-\log(|z|^2)$.
If $M$ is a complete K\"ahler manifold of finite volume and bounded negative 
sectional curvature, $M$ is hyper $1$--\,concave. This is shown by Siu--Yau in 
\cite{SY} by using Buseman functions. Moreover, if $\dim M\geqslant 3$, this example
falls in the previous case since by \cite{Nad} $M$ can be compactified to an algebraic space 
by adding finitely many points. 

\proclaim{Theorem 3.1}
Let $M$ be a hyper 1--\,concave manifold carrying a 
line bundle $(E,h)$ which is semi-positive outside a compact set. Let $\tm$ be a Galois covering
of group $\g$ and $\te$ the lifting of $E$. Then
$$\dg H^{n,0}_{(2)}(\tm, \te^k)\geqslant\frac{k^n}{n!}\int_{M(\leqslant 1,h)}
\left(\frac{\imath}{2\pi}{\boldkey c}(E, h)\right)^n +o(k^n)
\,,\quad k\longrightarrow\infty\,, $$
where the $L^2$ condition is with respect to $\widetilde h$ and any metric on $\tm$. 
\endproclaim

\demo{Proof}
Let us consider a proper function $\f :M\longrightarrow (-\infty,0\,)$ which is 
strictly plurisubharmonic outside a compact set. The fact that $\f$ goes to $-\infty$ 
to the ideal boundary of $M$ allows to construct a complete hermitian metric on 
$M$ which has moreover the feature of being K\"ahler outside a compact set. 
Namely we consider the function 
$\chi=-\log(-\f)$
so that
$
\pa\db\chi=\f^{-2}\,\pa\f\wedge\db\f-{\f}^{-1}\,\pa\db\f
$
which is obviously positive definite on the set where $\pa\db\f$ is. We can now patch 
$\pa\db\chi$ and an arbitrary hermitian metric on $M$ by using a smooth partition
of unity to get a metric $\omega_0$ on $M$ such that
$\omega_0=\pa\db\chi\quad\text{on}\quad M\setminus K,\;K\Subset M$.
It is easy to verify that $\omega_0$ is complete since
the function $-\chi$ is an exhaustion function and 
$
\omega_0=\omega+\pa(-\chi)\wedge\db(-\chi)
$
where $\omega=-\f^{-1}\pa\db\f$ is a metric on $M\setminus K$, so that $d(-\chi)$ is
bounded in the metric $\omega_0$\,. 
Note that $\omega_0$ is obviously K\"ahler on $M\setminus K$.

Let us consider a holomorphic hermitian line bundle $E$ endowed
with a metric $h$ such that $\imath{\boldkey c}(E,h)\geqslant 0$ on 
$M\setminus K$ (we stretch $K$ if necessary).
We equip $E$ with the metric $h_{\varepsilon}=h\exp(-\varepsilon\chi)$ and the
curvature relative to the new metric satisfies
$\imath{\boldkey c}(E,h_{\varepsilon})\geqslant\varepsilon\,\omega_0$
on $M\smallsetminus K$.
We are therefore in the conditions of Theorem 2.1. 
First observe that $h_\ve\gtrsim h$ so that $H_{(2)}^{n,0}(\tm,\te^k,\widetilde\omega_0,
\widetilde h_\ve)\subset 
H_{(2)}^{n,0}(\tm,\te^k,\widetilde\omega_0,\widetilde h)$ which is an injective $\g$--morphism. 
By Theorem 2.1
$$\liminf_{k}k^{-n}\dg H_{(2)}^{n,0}(\tm,\te^k,\widetilde\omega_0,\widetilde h_\ve)
\geqslant \frac{1}{n!}
\int_{\Omega(\leqslant 1,h_\varepsilon)}\left(\frac{\imath}{2\pi}{\boldkey c}(E, h_\ve)\right)^n 
$$ 
so that
$$
\liminf_{k}k^{-n}\dg H_{(2)}^{n,0}(\tm,\te^k)\geqslant \frac{1}{n!}
\int_{\Omega(\leqslant 1,h_\varepsilon)}\left(\frac{\imath}{2\pi}{\boldkey c}(E, h_\ve)\right)^n
\eqno(4.1)$$
We let now $\varepsilon\searrow 0$ in (4.1); since $h_\varepsilon$ converges uniformly together 
with its derivatives to $h$ on compact
sets we see that we can replace $h_\varepsilon$ with $h$
in the right-hand side of (4.1).
Let $M(q,h)$ be the set where $\imath{\boldkey c}(E, h)$ is non-degenerate and has exactly $q$ 
negative eigenvalues. By hypothesis $M(1,h)\subset K$ and on $M(0,h)=
M(\leqslant 1,h)\smallsetminus M(1,h)$ 
the integrand is positive. 
Hence we can let $\Om$ exhaust $X$ and we get the inequality from the statement of the theorem. 
\hfill{\,}\qed
\enddemo

\noindent
We prove now that Siu's criterion extends tale quale for hyper 1--concave manifolds.
\proclaim{Corollary 3.2}
Let $M$ be a hyper 1--\,concave manifold carrying a 
line bundle which is semi-positive outside a compact set and satisfies
Demailly's condition \rom{(D)}. Then $X$ is Moishezon.
In particular the conclusion holds true if $E$ is semi-positive and 
positive at one point. 
\endproclaim
\demo{Proof} By Theorem 3.1 (for $\g=\{\Id\}$) we have
$$\dim H^0(M,E^k\otimes K_M)\geqslant\dim H_{(2)}^{n,0}(M,E^k)
\geqslant C\,k^n$$ with $C>0$ for large $k$, by condition (D).
We note that the first space is finite dimensional since $M$ is $1$--concave.
By the Siegel--Serre Lemma (Proposition 5.7 from \cite{Ma}),
$\dim H^0(M,E^k\otimes K_M)\leqslant C\,{k^{\varkappa(E)}}$, ($k>0$),
where $\varkappa(E)$ is the supremum over $k$ of the generic 
rank of the canonical meromorphic mapping from
$M$ to $\Bbb P\big(H^0(M,E^k\otimes K_M)^*\big)$.
We obtain that
$\varkappa (E)=n$, that is, the line bundle $E^k\otimes K_M$ gives local
coordinates on an open dense set of $M$ for sufficiently large $k$.
This clearly implies $M$ Moishezon and thereby concludes the proof.
\hfill{\,}\qed
\enddemo
\remark{Remark 3.1}

\noindent
(a)
We can use this criterion in the Nadel 
compactification theorem \cite{Nad}. It asserts that if $M$ is a
connected manifold of dimension $\geqslant 3$ satisfying\,:
(i) $M$ is hyper 1--concave, (ii) $M$ is Moishezon,
(iii) $M$ can be covered by Zariski-open sets which are
uniformized by Stein manifolds,
then $M$ is biholomorphic to $M^*\smallsetminus S$ where $M^*$ is a compact
Moishezon space and $S$ is finite.
We see thus that condition (ii) in Nadel's theorem may be replaced with
the analytic condition: $M$ possesses a line bundle which is semi-positive
outside a compact set and satisfies Demailly's condition (D). 

\noindent
(b)
In general, if $M$ is a hyper $1$--\,convave
manifold of dimension $n\geqslant 3$ possesing a semi--positive line
bundle satistying (D) then (by a theorem of Rossi) it can be compactified
so that $M$ is biholomorphic to an open set of a compact
Moishezon manifold which is the complement of a complete pluripolar set. 
Therefore there exist a meromorphic mapping
defined on $X$ with values in a projective space which is an embedding
outside a proper analytic set of $X$. To see this we have to apply the
corresponding statement for compact Moishezon manifolds, a result due
to Moishezon. The difficulty in Nadel's theorem is to show that under
additional hypothesis the pluripolar set is actually a finite set.

\noindent
(c)
The argument in the proof of Corollary 3.2 shows that the 
integral appearing in Theorem 3.1
is finite. Thus, if $E$ is positive outside a compact set $K$ then $M\smallsetminus K$ has 
finite volume with respect to the 
metric $\imath{\boldkey c}(E)$ (this observation stems from \cite{NT}). 

\noindent
(d)
If $M$ possesses a positive line bundle $E$ then $\imath{\boldkey c}(E)+\imath\partial\db\chi$
is a complete K\"ahler metric and H\"ormander's $L^2$ estimates and Andreotti--Tomassini's theorem
\cite{AT} show that $E$ is ample and $M$ can be embedded in the projective space.
So even in dimension $2$ we can compactify $M$ (by \cite{An}).

\noindent
(f)
Let $X$ is a compact complex space of dimension $n\geqslant 2$ and with isolated
singularities. Suppose that we have a line bundle $E$ on $\rg$ which
is semi-positive in a deleted neighbourhood of $\sg$ and satisfies (D). Then $X$ is Moishezon. 
Indeed, by the previous result we find $n=\dim X$ independent meromorphic
functions on $\rg$ which extend to $X$ by the Levi extension theorem. 
This is a generalization of 
Takayama's criterion \cite{Ta} in the case of isolated singularities.
We allow weaker hypothesis, that is $E$ is defined just on $\rg$ and
the curvature condition is just semi-positivity. The reason is the
good exhaustion function we have at hand. In the general case one has
to use the Poincar\'e metric and the strict positivity near $\sg$ is
essential. Note however that the method of Takayama gives that the
line bundle who forms local coordinates is $E^k$, while in our proof
is $E^k\otimes K_X$.
\endremark

We want now to study the following type of stongly pseudoconcave manifold.
Let $X$ be an irreducible compact complex space with isolated singularities and 
of dimension $\geqslant 2$. We know that $\rg$ is hyper $1$--\,concave and we 
denote by $\varphi :\rg\longrightarrow\Bbb R$ the exhaustion function. Since 
$\varphi$ is strictly plurisubharmonic outside a compact set we have that the sub--level
sets $X_c=\{\varphi >c\}$ are $1$--\,concave manifolds 
i.e. stongly pseudoconcave domains. 
In our previous paper \cite{Ma} 
we have shown that in general if $M$ is a $1$--concave manifold of dimension $\geqslant 3$
which carries a hermitian line bundle $E$ which semi-negative near the boundary and satisfies 
(D) then the Kodaira dimension of $E$ is maximal and $M$ is Moishezon.
The assumption about the change of curvature sign
(i.e. semi-negativity) near the boundary is imposed
by the construction of complete hermitian metrics $\omega_0$ and $h_0$ as in Theorem 2.2
which give the $L^2$ estimate and preserve condition (D) for $h_0$; the negativity of the Levi form
of the sublevel sets of $M$ requires as a natural curvature condition for $E$
the semi--negativity. The restriction on dimension comes from the fact that we need an $L^2$
estimate in bi--degree $(0,1)$. 
Of course, usually we are given an overall positive bundle $E$ on $M$. 
We show that for manifolds $X_c$ as before we can also apply the
criterion in \cite{Ma} alluded to by modifying the metric. 

We recall at the outset some terminology.
Let us consider a covering $\{U_\alpha\}$ of $X$ and embeddings
$\iota_\alpha:U_\alpha\hookrightarrow\Bbb C^{N_\alpha}$ such that
$E|_{U_\alpha}$ is the inverse image by $\iota_\alpha$ of the trivial
line bundle $\Bbb C_\alpha$ on $\Bbb C^{N_\alpha}$. 
Moreover we consider hermitian metrics $h_\alpha=e^{-\varphi_\alpha}$ 
on $\Bbb C_\alpha$ such that $\iota_{\alpha}^*{h_\alpha}=
\iota_{\beta}^*{h_\beta}$ on $U_\alpha\cap U_\beta\cap\rg$. The system
$h=\{\iota_{\alpha}^*{h_\alpha}\}$ is called a hermitian metric on $E$ over $X$. It
clearly induces a hermitian metric on $E$ over $\rg$. 
\comment
We shall allow our
metrics to be singular at the singular points, that is, $\f_\alpha\in
L^1_{loc}(\Bbb C^{N_\alpha})$ is smooth outside
$\iota_\alpha(\sg)$. 
\endcomment
The curvature current $\imath{\boldkey c}(E)$ is given in 
$U_\alpha$ by $\iota_{\alpha}^*(\imath{\pa\db\f_\alpha})$ which on $\rg$
agrees with the curvature of the induced metric.

\comment
We shall suppose in
the sequel that the curvature current is dominated by the euclidian metric i.e. 
$\imath\pa\db\f_\alpha$ is bounded above and below by constant
times $\omega_E=\imath\sum dz_j\wedge d{\bar z}_j$.
\endcomment

\proclaim{Theorem 3.3}
Let $X$ be an irreducible compact complex space with isolated singularities and let
$X_c$ be the sublevel sets of the hyper $1$--\,concave manifold $\rg$. Assume that there 
exists a holomorphic line bundle $E\longrightarrow X$ with a smooth hermitian metric such that
condition \rom{(D)} is fulfilled on $\rg$. Then for sufficiently small $c$ there exists a metric on $E$
such that $E$ is negative in the neighbourhood of $bX_c$ and $\int_{X_c(\leqslant 1)}
\big(\imath{\boldkey c}(E)\big)^n>0$.
\endproclaim

\demo{Proof}
Let $\pi :\widetilde X\longrightarrow X$ be a resolution of
singularities of $X$. Let us denote by $D_i$ the components of the
exceptional divisor. Then there exist positive integers $n_i$ such
that $D:=\sum n_i\,D_i$ admits a smooth hermitian metric such that
the induced line bundle $[D]$ is negative in a neighbourhood
$\widetilde U$ of $D$ (cf. \cite{Sa}). Let us consider a canonical
section $s$ of $[D]$, i.e. $D=(s)$, and denote by $|s|^2$ the pointwise
norm of $s$ with respect to the above metric. By Lelong-Poincar\'e equation
$\f=\log|s|^2$ is strictly plurisubharmonic on $\widetilde
U\setminus D$. By using a smooth function on $\widetilde X$ with 
compact support in $\widetilde U$ which equals one near $D$ we
construct a smooth function $\chi$ on $\widetilde X\smallsetminus D\simeq\rg$ such that 
$\chi=-\log(-\log|s|^2)$ on $\widetilde U\setminus D$.   

\noindent
Since $\log|s|^2$ goes to $-\infty$ on $D$, this is the analogue of the
function constructed in the proof of Theorem 3.1\,. As there we show
that $\imath\pa\chi\wedge\db\chi\leqslant\imath\pa\db\chi$. 
Let us consider a metric $\omega$ on $\rg$ which on every open set
$U_\alpha$ as above is the pullback of a hermitian metric on the
ambient space $\Bbb C^{N_\alpha}$,
$\omega=\iota_{\alpha}^*\,{\omega_\alpha}$\,.
We consider then the metric (K\"ahler near $\sg$)
$\omega_0=A\omega+\pa\db\chi$
where $A>0$ is chosen sufficiently large (to ensure that $\omega_0$ is
a metric away from the open set where $\pa\db\chi$ is positive definite).
It is easily seen that $\omega_0$ is complete by the same argument as
in the proof of Theorem 3.1\,. This kind of metrics were introduced by
Saper in \cite {Sa}. They have finite volume.

Let us consider now a neighbourhood $U$ of the singular set. We
assume that $U$ is small enough so that there are well defined on $U$
a potential $\rho$ for $\omega$ and a potential $\eta$ for the
curvature $\imath{\boldkey c}(E)$ (they are restrictions from ambient spaces). 
By suitably cutting-off we may define a function $\psi\in{\Cal C}^\infty
(\rg)$ such that
$\psi=-\chi-\eta-A\,\rho$
near $\sg$\,. Remark that since $\imath{\boldkey c}(E)$ is bounded
above by a continuous $(1,1)$ form near $\sg$ the potential $-\eta$ is
bounded above near the singular set. This holds true for $\rho$ too
(it is smooth) so that $\psi$ tends to $\infty$ at the singular 
set $\sg$. Let us consider a smooth function $\gamma :\Bbb
R\longrightarrow\Bbb R$ such that 
$$\gamma (t)=\cases
0&\text{if} \quad t\leqslant 0\,,\\
t&\text{if}\quad t\geqslant 1\,.
\endcases$$
and the functions  $\gn :\Bbb R\longrightarrow\Bbb R$  
given by $\gn(t)=\gamma(t-\nu)$ for all positive integers $\nu$\,.
Let us denote the hermitian metric on $E$ by $h$ and let us consider
the following metric on $E$\,: 
$h_\nu=h\,\exp\big(-\,\gn(\psi)\big)$,
with curvature
$$\imath{\boldkey c}(E,h_\nu)=\imath{\boldkey
c}(E,h)+\gpn(\psi)\pa\db\psi + \gsn(\psi)\pa\psi\wedge\db\psi\,.$$
On the set $\{\psi\geqslant\nu+1\}$ we have $\gn(\psi)=
\psi-\nu$ so that $\gpn(\psi)=1$ and $\gsn(\psi)=0$ and therefore
$\imath{\boldkey c}(E,h_\nu)=\imath{\boldkey c}(E,h)+\pa\db\psi$.
Since $\psi$ goes to $\infty$ when we approach the singular set we
may choose $\nu_0$ such that for $\nu\geqslant\nu_0$ we have 
$\{\psi\geqslant\nu+1\}\subset U$ where $U$ is a sufficiently small neighbourhood of
$\sg$. Bearing in mind the meaning of
$\eta$ and $\rho$ together with the definition of $\omega_0$ it is
straightforward that
$\imath{\boldkey c}(E,h_\nu)=-\omega_0$ on 
$\{\psi\leqslant\nu+1\}$, that is $(E,h_\nu)$ is negative on this set.
We denote $\Omega_\nu$ the compact set $\{\psi\leqslant\nu+2\}$\,.
We decompose this set in $\Omega^\prime_\nu=\{\psi\leqslant\nu\}$ and
$\Omega^{\prime\prime}_\nu=\{\nu\leqslant\psi\leqslant\nu+2\}$. 
On $\Omega^\prime_\nu$ we have $\gn(\psi)=0$ and 
$\imath{\boldkey c}(E,h_\nu)=\imath{\boldkey
c}(E,h)$\,. We infer that
$$\eqalignno{
\int_{\Omega^\prime_\nu(\leqslant 1,h_\nu)}\big(\imath{\boldkey c}
(E,h_\nu)\big)^n\,&=\,\int_{\Omega^\prime_\nu
(\leqslant 1,h)}\big(\imath{\boldkey c}(E,h)\big)^n\cr
&=\int_{\rg(\leqslant 1,h)}{\boldkey
1}_{\Omega^\prime_\nu}\alpha_1\dotsm\alpha_n\,dV_0&(4.2)}$$
where $\alpha_1,\dotsc,\alpha_n$ are the eigenvalues of 
$\imath{\boldkey c}(E,h)$ with respect to $\omega_0$ and $dV_0$ is the
volume form of the same metric. Since 
$\imath{\boldkey c}(E,h)$ is dominated by the euclidian metric near $\sg$, 
$\imath{\boldkey c}(E,h)$ is dominated by $\omega$ and by $\omega_0$.
Hence the product $\alpha_1\dotsm\alpha_n$ is bounded on $\rg$. Since 
$\rg(\leqslant 1)$ has finite volume with respect to $\omega_0$ the
functions $|{\boldkey
1}_{\Omega^\prime_\nu}\alpha_1\dotsm\alpha_n|$ are bounded by an
integrable function. On the
other hand ${\boldkey 1}_{\Omega^\prime_\nu}\longrightarrow 1$ when 
$\nu\longrightarrow\infty$ so that the integrals in (4.2) tend to 
$\int_{\rg(\leqslant 1,h)}\big(\imath{\boldkey c}
(E,h)\big)^n$ which is assumed to be positive. 
Thus it suffices to show that the integral on the set
$\Omega^{\prime\prime}_\nu$\;\, i.e. $\int_{\Omega^{\prime\prime}_\nu
(\leqslant 1,h_\nu)}\big(\imath{\boldkey c}(E,h_\nu)\big)^n$
tends to zero as $\nu\longrightarrow\infty$. For this purpose we use
the obvious bound
$$\int_{\Omega^{\prime\prime}_\nu
(\leqslant 1,h_\nu)}\left(\frac{\imath}{2\pi}{\boldkey c}(E,h_\nu)\right)^n\,
\leqslant\sup |\,\delta_1\dotsm\,\delta_n|\cdot\operatorname{vol}\, 
(\Omega^{\prime\prime}_\nu)$$
where $\delta_1,\dotsc,\delta_n$ are the eigenvalues of
$\imath{\boldkey c}(E,h_\nu)$ with respect to $\omega_0$ and the
volume is taken in the same metric. We use now the minimum-maximum
principle to see that: 
(i) $\delta_1$ is bounded below and
$\delta_2,\dotsc,\delta_n$ are bounded above on the set of integration 
$\Omega^{\prime\prime}_\nu(1,h_\nu)$ and (ii) $\delta_1,\dotsc,\delta_n$
are upper bounded on $\Omega^{\prime\prime}_\nu(0,h_\nu)$.
For this we need the
domination of $\imath{\boldkey c}(E,h)$ by $\omega$ and the
boundedness of $\gamma^\prime_\nu$ and $\gamma^{\prime\prime}_\nu$\,.
Since $\operatorname{vol}\,(\Omega^{\prime\prime}_\nu)\longrightarrow
0$ as $\nu\longrightarrow\infty$ our contention follows.
Hence for large $\nu$ the metric $h_\nu$ does the required job. 
\hfill{\,}\qed 
\enddemo

\remark{Remark 3.2}
We have seen that Siu's criterion generalizes to compact complex
spaces with isolated singularities. Demailly's
criterion extends too. 
Let $X$ be an irreducible compact complex space of dimension
$n\geqslant 3$ with isolated
singularities and $E$ a smooth hermitian line bundle over $X$.
Assume that condition \rom{(D)} is fulfilled on $\rg$. Then $X$ is Moishezon.
Indeed, for small $c$ the sets $X_c$ are Moishezon by Corollary 4.3 of \cite{Ma}
and the meromorphic functions from $X_c$ extend to $X$.
In fact the result holds also for $n=2$ with a proof very similar to that of Theorem 4.4.
We note also that we can allow the metric $h$ of $E$ to be singular at $\sg$ but the
cuvature current $\imath{\boldkey c}(E)$ should be dominated (above ans below) by the euclidian
metric near $\sg$. The proof of Theorem 4.4 goes through with minor changes.
\endremark

Since the manifold $\overline{X}_c$ is compact Theorem 4.4 can be used to prove some 
stability results for certain perturbation of the complex structure of ${\overline X}_c$.
Since our approach relies on the use of a sufficiently positive line bundle $E$ we need 
to consider perturbations of the complex structure which lift to a perturbation of $E$.
This kind of sitiuation was studied by L. Lempert in \cite{Le}.
\proclaim{Proposition 3.4} Let $X$ be a Moishezon variety with isolated singularities and dimension
$n\geqslant 3$.
Let $\Cal J$ denote the complex structure of $\rg$ and
let $Z\subset\rg$ be a non--singular hypersurface such that the line bundle $E=[Z]$ satisfies \rom{(D)}.
Then for sufficiently small $c$ and any complex structure $\Cal J^\prime$ on ${\overline X}_c$
such that
\roster
\item $T(Z)$ is $\Cal J^\prime$ invariant and  
\item $\Cal J^\prime$ is sufficiently close to $\Cal J$ in the $\Cal C^\infty$ topology
\endroster
there exists a $\Cal J^\prime$--holomorphic line bundle $E^\prime$ on ${\overline X}_c$ 
which is negative near $b X_c$ and satisfies \rom{(D)}.
In particular $({\overline X}_c,{\Cal J^\prime})$ is a Moishezon pseudoconcave manifold 
and any compactification of $({\overline X}_c,{\Cal J^\prime})$ is Moishezon.
\endproclaim

\demo{Proof}
Let us first choose $c_0$ such that for $c<c_0$ there exists a `good' hermitian metric $h$ on 
$E$ over a neighbourhood of $X_c$, that is, whith negative curvature near the boundary 
and satisfying (D). We use now the description of the lifting of $\Cal J^\prime$ with properties
(1) and (2) as given in \cite{Le}. Namely, $Z$ determines a new $\Cal J^\prime$ holomorphic line bundle
$E^\prime\longrightarrow ({\overline X}_c,{\Cal J^\prime})$. There exists a finite open covering
${\Cal U}=\{U\}$ of ${\overline X}_c$ such that $E$ and $E^\prime$ are trivial on each $U$
and they are defined by multiplicative cocycles $\{g_{UV}\,{\Cal J}\,\text{holomorphic on}\,
{\overline U}\cap{\overline V}\;:\;U,V\in\Cal U\}$ and 
$\{g^\prime_{UV}\,{\Cal J^\prime}\,\text{holomorphic on}\,{\overline U}\cap{\overline V}\;:\;U,V\in\Cal U\}$.
Moreover $g_{UV}$ and $g^\prime_{UV}$ are as close as we please assuming $\Cal J$ and $\Cal J^\prime$
are sufficiently close. (By `close' we always understand close in the $\Cal C^\infty$ topology.)
Next we can define a smooth bundle isomorphism $E\longrightarrow E^\prime$ by resolving the smooth 
additive cocycle $\log(g^\prime_{UV}/g_{UV})$ in order to find smooth functions $f_U$, close to 1 
on a neighbourhood of $\overline U$ such that $g^\prime_{UV}=f_U\,g_{UV}\,f^{-1}_V$.
Then the isomorphism between $E$ and $E^\prime$ is given by $f=\{f_U\}$. 
The metric $h$ is given in terms of the covering $\Cal U$ by a collection $h=\{h_U\}$ of smooth
strictly positive functions satisfying the relation $h_V=h_U\,|g_{UV}|$. We define a hermitian metric
$h^\prime=\{h^\prime_U\}$ on $E^\prime$ by $h^\prime_U=h_U\,|f^{-1}_U|$; $h^\prime_U$ is close to $h_U$.
The curvatures forms of $E$ and $E^\prime$ are given by 
$$\frac{\imath}{2\pi}{\bold c}(E)=\frac{1}{4\pi}\,d\circ{\Cal J}\circ d\,(\log h_U)\,,\quad
\frac{\imath}{2\pi}{\bold c}(E^\prime)=\frac{1}{4\pi}\,d\circ{\Cal J^\prime}\circ d\,(\log h^\prime_U)\,.$$
Therefore, when $\Cal J^\prime$ is sufficiently close to $\Cal J$, $\frac{\imath}{2\pi}{\bold c}(E^\prime)$
is negative near the boundary of ${\overline X}_c$ and, since the eigenvalues of 
$\frac{\imath}{2\pi}{\bold c}(E^\prime)$  are close to those of $\frac{\imath}{2\pi}{\bold c}(E)$,
$E^\prime$ satisfy the condition (D) i.e.  $\int_{X_c(\leqslant 1)}\big(\imath{\boldkey c}(E^\prime)\big)^n>0$. 
We can apply thus the Corrolary 4.3 of \cite{Ma} to the strongly pseudoconcave manifold
$({\overline X}_c,{\Cal J^\prime})$ to conclude that $({\overline X}_c,{\Cal J^\prime})$ is Moishezon.
\hfill{\,}\qed
\enddemo

\remark{Remark 3.3}
If $[Z]$ is positive, part of the stability property follows from the rigidity of embeddings with 
positive normal bundle. Indeed, assume $N_{Z}=[Z]\upharpoonright_Z$ 
is positive in $({\overline X}_c,{\Cal J^\prime})$ 
(for any $c$ such that this manifold is still pseudoconcave). 
Then Ph.~Griffiths \cite{Gri1} has shown that 
there exists a neighbourhood $W$ of $Z$ such that the mapping $\Phi:(X_c,{\Cal J}^\prime)- - \rightarrow{\Bbb P}^N$
given by $[mZ]$ is an embedding of $W$ for large $m$ . Thus $(X_c,{\Cal J}^\prime)$ is Moishezon.
Our result deals with the slightly more general situation of a `big' embedding
i.e. when $[Z]$ is not ample but satisfies condition (D). Moreover we have a useful quantitative way
of measuring whether the perturbed structure is Moishezon.
\endremark
\proclaim{Corollary 3.5}
Let $({\overline X}_c,\Cal J^\prime)$ and $E^\prime$ be as in Proposition 4.6. 
Then there exists hermitian metrics on $X_c$ and $E^\prime$ and a positive constant $C$ such that 
for any Galois covering ${\widetilde X_c}\longrightarrow X_c$ of group $\g$ we have 
$$\dim_{\sG}H^0_{(2)}({\widetilde X_c}\,,\widetilde{E^\prime}^k)\geqslant C\,k^n+o(k^n)
\,,\quad k\longrightarrow\infty\,.$$ 
the $L^2$ condition being with respect to lifts of the hermitian metrics on $X_c$ and $E^\prime$.
\endproclaim
\demo{Proof}We know that we have on $E^\prime$ a metric $h$ satisfying the conclusion
of Theorem 4.4. 
Then, as in Theorem 2.2, we can construct metrics $\omega_0$ and $h_0$ in order to obtain (2.9).
Note that the integral in (2.9) depends on the modified metric $h_0$ so we cannot always infer that 
it is positive even if $(E^\prime,h)$ satisfies (D). But under the assumption of semi--negativity of $h$
near the boundary we can construct an $h_0$ such that the integral in (2.9) is positive (cf. 
Corollary 4.3 of \cite{Ma}). Thus by applying Theorem 2.2 we get the conclusion.
\hfill{\,}\qed
\enddemo

\subhead
\S 4 $L^2$ generalization of a theorem of Takayama 
\endsubhead

In this section we study
the $L^2$ cohomology of coverings of Zariski open sets in compact complex spaces.
For compact spaces with singularities 
Takayama \cite{Ta} generalized Siu--Demailly criterion if 
$E\longrightarrow X$ is a line bundle endowed with a singular
hermitian metric which is smooth outside a proper analytic set
$Z\supset\sg$ and defines a strictly positive current near $Z$. 

Using the setting of Takayama's theorem we shall study coverings of Zariski open sets
in compact complex spaces. 
\proclaim{Proposition 4.1}
Let $X$ be an $n$--dimensional compact manifold
and let $E$ be a holomorphic line bundle with a singular hermitian metric $h$.
We assume that:
\roster
\item      $\imath{\boldkey c}(E,h)$ is smooth on $M=X\smallsetminus Z$ where 
$Z$ is a divizor with only simple normal crossings\,;
\item      $\imath{\boldkey c}(E,h)$ is a strictly positive current in a neighbourhood
of $Z$\,.
\endroster
\noindent
Let $p:\tm\longrightarrow M$ be a Galois covering with group
$\g$ and $\te=p^\ast E$. 
Then, 
$$\dg H^0_{(2)}(\tm,\te^k)\geqslant \frac{k^n}{n!}\int_{M(\leqslant 1,h)}
\left(\frac{\imath}{2\pi}{\boldkey c}(E, h)\right)^n +o(k^n)
\,,\quad k>>0\,, $$
where $H^0_{(2)}(\tm,\te^k)$ is the space of sections 
of $\te^k$ which are $L^2$ with respect to the pullbacks of the restrictions to $M$ and
$E\upharpoonright _M$ of  
smooth metrics on $X$ and $E$.
\endproclaim

\demo{Proof}
This is an equivariant form of Takayama's main technical result in \cite{Ta}.
Namely we construct the Poincar\'e metric $\omega_\ve$ on $M$ (for details see \cite{Zu})
and $h_\ve$ as in \cite{Ta} and remark that the hypothesis of Theorem 2.1 are satisfied.
Moreover we can work with $(0,1)$--forms since the Ricci curvature of the Poincar\'e metric
is bounded below.

More specifically, we write $Z=\sum Z_j$ and consider a section $\sigma_j$ of the line bundle
$[Z_j]$ which vanishes to first order on $Z_j$. Then we endow $[Z_j]$ with 
a hermitian metric such that the norm of $\sigma_j$ satisfies $|\sigma_j|<1$.
Take then an arbitrary smooth metric $\omega^\prime$ on $X$ and define
$\omega_\ve =\omega^{\prime}-\ve\,\imath\,\sum\partial\db
(-\log|\sigma_j|^2)^2$ on $M=X\smallsetminus Z$
which for small $\ve>0$ is a complete metric on $M$.
Then we consider the following family of metrics on $E\upharpoonright _M$:
$h_\ve =h\prod_j (-\log |\sigma_j|^{2})^{\ve}$, $\ve >0$.
We check now the hypotheses of Theorem 2.1 is satisfied. First we remark that
the torsion operators of the
Poincar\'e metric are pointwise bounded with respect to the Poincar\'e metric since $d\omega_\ve=
d\omega^\prime$ and $\omega_\ve > \omega^\prime$. Also the Ricci curvature ${\boldkey c}(K^{\ast}_{\tm})$
of $\widetilde\omega_\ve$ is bounded below with respect to $\widetilde\omega_\ve$ by a constant 
independent of $\ve$ (since this is true for $\omega_\ve$). 
Since $E$ is strictly positive in the neighbourhood of $Z$ condition (A) is satisfied for a compact
$K$ outside which $E$ is positive (and it doesn't depend on $\ve$).

Let $h^\prime$ be a smooth hermitian metric on $E$ over $X$. Near $Z$ the metric $h$
is locally represented by a strictly plurisubharmonic weight. Thus $h$ is locally
bounded below near $Z$ and thus $h\geqslant C\,h^\prime$ on $X$ for some positive 
constant $C$.
We remark now that $h_\ve >h\geqslant C\,h^\prime$ and $\omega_\ve>\omega^\prime$ 
near $Z$ so that we have the inclusion 
$H^0_{(2)}(\tm,\te^k)_\ve\subset H^0_{(2)}(\tm,\te^k)$, 
(which is an injective $\g$--morphism)
in the last group the $L^2$ condition being taken with respect to $\widetilde h^\prime$
and $\widetilde\omega^\prime$. 
By Theorem 2.1 for $K\Subset\Om\Subset M$
$$
\dg H^0_{(2)}(\tm,\te^k)\geqslant\dg H^0_{(2)}(\tm,\te^k)_\ve 
\geqslant\int_{\Om (\leqslant 1,h_{\ve})}
\left(\frac{\imath}{2\pi}{\boldkey c}(E, h_\ve)\right)^n
+o(k^n)\,.$$
We can let $\ve\longrightarrow 0$ in the right--hand side in order to replace $h_\ve$ with $h$.
Then we can let $\Om$ exhaust $X$ to get the inequality from the statement. 
\hfill{\,}\qed
\enddemo

\proclaim{Theorem 4.2}
Let $X$ be an irreducible reduced compact Moishezon space and let $M\subset
{\operatorname{Reg}}(X)$ be a
Zariski open set. There exists a holomorphic line bundle $E\longrightarrow\rg$
endowed with a singular hermitian metric whose curvature current $\imath{\boldkey c}(E)$
is positive and such that for any Galois covering $\tm @>p>>M$ of group $\g$ we have
$$
\dg H^0_{(2)}(\tm,\te^k)\geqslant\frac{k^n}{n!}\int_M \left(\frac{\imath}
{2\pi}{\boldkey c}(E)\right)^n +o(k^n)\,,\quad k\longrightarrow\infty
$$ 
where the integration takes place outside $\operatorname{Sing\,supp}{\boldkey c}(E)$.
The $L^2$ condition is taken with respect to liftings of smooth hermitian 
metrics on $M$ and $E$ induced from a resolution of singularities of $X$.
\endproclaim
\demo{Proof}
{\it Step 1.}\;Let $X$ be a Moishezon manifold and $M=X\setminus Z$ a Zariski open set, where
$Z$ is a proper analytic set. Thanks to Moishezon $X$ admits a projective modification.
Therefore there exists a strictly positive integral K\"ahler
current $T$ on $X$. Equivalently there exists a holomorphic line bundle $E$ on $X$
possesing a singular hermitian metric such that the curvature current $T=\imath{\boldkey c}(E)$
is strictly positive (bounded below by a smooth hermitian metric). Assume that
$\operatorname{Sing\,supp}T\subset Z$. Then $M$ is biholomorphic to a Zariski open set as in the 
statement of Proposition 3.1. Indeed, we can blow up $Z$ to make it a divisor with only
simple normal crossings. By replacing $E$ with higher tensor powers and twisting it with
the dual of the exceptional divisor at each step of the blowing up process we can ensure
that on the blow--up we still have a positive line bundle with singular metric along $Z$.
Thus in this case we can apply Proposition 3.1.

\noindent
{\it Step 2.}\;
To go further
let $M$ be a Zariski open set in a Moishezon manifold $X$.
By a theorem of Demailly \cite{De2} we know that there exists a strictly positive integral 
K\"ahler current $T$ with analytic singularities. As a consequence 
$\operatorname{Sing\,supp}T\subset S$, where $S$ is a proper analytic set. As before we 
can suppose that $S\cup Z$ is a divisor with only simple normal crossings. Let $E$ be a 
line bundle with singular hermitian metric such that $T=\imath{\boldkey c}(E)$.
Denote by $M_1=X\setminus (S\cup Z)=M\setminus S$\,: $M_1$ and $E$ are as in Proposition 3.1. 
Let $p:\tm\longrightarrow M$ be a Galois covering of group $\g$. Setting
$\tm_1=p^{-1}M_1$ we have a Galois covering $\tm_1\longrightarrow M_1$ of group $\g$.
Hence, $\dg H^0_{(2)}(\tm_1,\te^k)\geqslant\frac{k^n}{n!}\int_{M_1} \left(\frac{\imath}
{2\pi}{\boldkey c}(E)\right)^n +o(k^n)$, for $k\longrightarrow\infty$.
The $L^2$ condition on $\tm_1$ is with respect to liftings of smooth 
hermitian metrics on $X$ and $E$. But a holomorphic section defined outside the 
analytic set $\widetilde S=p^{-1}S$ which is square integrable with respect 
to a smooth metric on 
$\tm$ extends past $\widetilde S$ as a holomorphic section on $\tm$.
We infer $\dg H^0_{(2)}(\tm,\te^k)\geqslant\frac{k^n}{n!}\int_{M} \left(\frac{\imath}
{2\pi}{\boldkey c}(E)\right)^n +o(k^n)$
the integral being taken on the smooth locus of $\imath{\boldkey c}(E)$
The $L^2$ condition is taken with respect to pullbacks of smooth metrics
on $X$ and $E$. 

\noindent
{\it Step 3.}\;
Finally let $X$ and $M$ as in hypothesis. By a resolution of singularities $M$ is biholomorphic
to a Zariski open set of a Moishezon manifold. By the preceding remarks we can conclude.
\hfill{\,}\qed
\enddemo

The following Proposition is a consequence of Theorem 4.2 in the case of Galois coverings
(taking into acount that the number of sheets of such a covering equals the cardinal of $\g$).
However, using Theorem 2.2 of Napier and Ramachandran \cite{NR}, we
can prove it for any unramified covering.
\proclaim{Proposition 4.3}
Let $X$ be an irreducible reduced compact Moishezon space and let $M\subset
{\operatorname{Reg}}(X)$ be a
Zariski open set. There exists a holomorphic line bundle $E\longrightarrow\rg$
such that for any unramified covering $p:\tm\longrightarrow M$ we have
$$\dim H^0_{(2)}(\tm,\te^k)\geqslant C\,k^n\,d\,,\quad k>>0\eqno{(4.1)}$$    
where $d$ is the number of sheets of the covering and $C>0$.
\endproclaim
\demo{Proof}
In the situation of Step 1 of the preceding proof 
 we see that the Poincar\'e metric on
$M$ is a complete K\"ahler metric since $M$ has the K\"ahler metric $\imath{\bold c}(E)$.
Therefore $\tm$ possesses a complete K\"ahler metric and a positive line bundle $\te$.
By applying the $L^2$ estimates of H\"ormander as in \cite{NR, Theorem 2.2} we get the result
for the $L^2$ cohomology with respect to the metrics $\widetilde\omega_\ve$ and $\widetilde h_\ve$
(notations of Proposition 4.1).
As in the proof of Proposition 4.1 we see that we can actually use pull-backs of smooth metrics on $X$.
Steps 2 and 3 go through as before.
\hfill{\,}\qed
\enddemo
\subhead
\S 5 Further remarks
\endsubhead
\comment
We can apply Theorem 2.1 to the case of a complete K\"ahler manifold $M$ 
of bounded negative Ricci curvature; see Nadel and Tsuji \cite{NT}. Let us remark that by a theorem of
Mok and Yau there exists a unique K\"ahler--Einstein metric, invariant by analytic automorphisms,
on any bounded domain of holomorphy in $\Bbb C^n$. Thus this metric descends to a complete
K\"ahler--Einstein metric on any quotient of the domain by a properly discontinuous discrete group.
\proclaim{Proposition 5.1}
Let $(M,\omega)$ be a complete K\"ahler manifold such that $\operatorname{Ric}\omega$
is negative and bounded above outside a compact set. Then
$$\dg H^0_2(\tm, K^k_{\tm})\geqslant \frac{k^n}{n!}\int_{M(\leqslant 1)}
\left(-\frac{\imath}{2\pi}\operatorname{Ric}\omega\right)^n +o(k^n)
\,,\quad k\longrightarrow\infty\,,$$ 
for any Galois covering $\tm$ of $M$. The $L^2$ condition is with respect
to the lift of the meric $\omega$ and $M(\leqslant 1)$ is the set where
$-\operatorname{Ric}\omega$ is non--degenerate and has at most one negative eigenvalue.
\endproclaim
We have ${\boldkey c}(K_M)=-\operatorname{Ric}\omega$ so a straightforward application of 
Theorem 2.1 gives the conclusion. Note that the case $\g=\{\Id\}$ and 
$\operatorname{Ric}\omega$
everywhere negative bounded above has been proved by Nadel and Tsuji.
\proclaim{Corollary 5.2}
Assume $M$ admits a complete K\"ahler--Einstein metric $\omega$ with 
$\operatorname{Ric}\omega=-\omega$ and of infinite volume.
Then $\dg H^0_2(\tm, K^k_{\tm})=\infty$ for $k$ large enough.
\endproclaim 
We remark that the last conclusion is stronger than the results coming from
the $L^2$ method which give just $\dim H^0_2(\tm, K^k_{\tm})\geqslant C\,|\g|\,k^n$
for some positive constant $C\in\Bbb R$.
\endcomment

We will apply Theorem 2.1 to the case of a complete K\"ahler manifold $M$ with positive
canonical bundle $K_M$. The case $\g=\{\Id\}$ is due to Nadel and Tsuji \cite{NT}.
If $D$ is a bounded domain of holomorphy in $\Bbb C^n$ we know by a theorem of Bremermann
that the Bergman metric $\omega=\omega_B$ is complete. On the other hand the Bergman metric is invariant under analytic automorphisms. Thus this metric descends to a complete
K\"ahler metric on any quotient of the domain by a properly discontinuous discrete group
$\g\subset\operatorname{Aut}(D)$. We denote $M=D/\g$ and
$\omega_\ast$ the induced Bergman metric on $M=D/\g$. If we denote by $B(z,\overline z)$
the Bergman kernel of $D$ we know that $B^{-1}$ can be considered as a hermitian
$\g$--invariant metric on $K_D$. Since $\omega=\partial\db\log B(z,\overline z)$ there
exists a hermitian metric on $K_M$ such that ${\boldkey c}(K_M)=\omega_\ast$. 
We have thus the following.
\proclaim{Proposition 5.1}
Let  $D$ is a bounded domain of holomorphy in $\Bbb C^n$, $\g\subset\operatorname{Aut}D$ a discrete 
group acting properly discontinuously on $D$ and $M=D/\g$. Then
$$\dg H^0_2(D, K^k_{D})\geqslant \left(\frac{k}{2\pi}\right)^n\int_{M}
\frac{\omega_\ast^n}{n!} +o(k^n)\,,\quad k\longrightarrow\infty\,$$
where the $L^2$ condition is taken with respect to the Bergman metric on $D$ and the metric $B^{-1}$
on $K_D$.
\endproclaim
Note that the space $H^0_2(D, K^k_{D})$ is a space of square integrable functions with respect to
the Bergman metric and to the weight $B^{-k}$.
An immediate consequence is the following.
\proclaim{Corollary 5.2} 
Assume that the Bergman metric on $M$ has infinite volume.
Then $\dg H^0_2(D, K^k_{D})=\infty$ for $k$ large enough.
\endproclaim

\noindent
We remark that the last conclusion is stronger than the results coming from
the $L^2$ method which gives just $\dim H^0_2(D, K^k_{D})\geqslant C\,|\g|\,k^n$
for some positive constant $C\in\Bbb R$.

Let us see what become our results in the simplest case of the unit disk $D\subset\Bbb C$.
Then the Bergman metric equals the hyperbolic metric $(1-|z|^2)^{-2}dz\wedge d\overline z$.
If $\g$ is a Fuchsian group, we have the following possibilities for large $k$:
\item{(a)} If $M=D/\g$ is compact, $\dg H^0_2(D, K^k_{D})= k\operatorname{vol}(M) +o(k)$.
\item{(b)} If $M$ is non--compact and has a finite number of cusps, 
the hyperbolic volume $\operatorname{vol}(M)$ is finite and
$\dg H^0_2(D, K^k_{D})\geqslant k\operatorname{vol}(M)+o(k)$,
\item{(c)} If $M$ is non--compact and the discontinuity set $\Omega\subset S^1$
is a union of intervals, $\dg H^0_2(D, K^k_{D})=\infty$ (since $\operatorname{vol}(M)=\infty$).
  
According to a conjecture of 
Griffiths \cite{Gri2, p.50}, if $D$ is a bounded domain in $\Bbb C^n$ which is topologically a cell
and $D/\g$ is quasi--projective then (i) the Bergman metric on $D/\g$ is complete and (ii)
the volume of $D/\g$ with respect to this metric is finite. In the sequel we discuss
the conjecture without the topological restriction.
If $D$ is a domain of holomorphy and $M=D/\g$ is pseudoconcave 
(e.g. $\operatorname{codim}(\overline 
M\smallsetminus M)\geqslant 2$), the answer is yes. Indeed, this follows from the Riemann--Roch 
inequalities for $\g=\{\Id\}$ in Proposition 5.1.
If $D$ is not necessarily a domain of holomorphy but $D/\g$ can be compactified
by adding a finite number of points we can show that the answer to (ii) is affirmative.
We do not assume $D/\g$ quasi--projective.
\proclaim{Proposition 5.3}
Let $D\Subset\Bbb C^n$ be an open set having a properly discontinous group $\g\subset\operatorname
{Aut}D$ such that there exists a compact complex space $Y$ with $D/\g\subset\operatorname{Reg}Y$ 
and $D/\g=Y\smallsetminus S$, where $S$ is a finite set. Then the volume of $D/\g$ in the 
induced Bergman metric is finite. 
\endproclaim
\demo{Proof}
Since $M=D/\g$ is hyper 1--concave and possesses a positive canonical bundle, we may apply 
Theorem 3.1 for $\g$ trivial and $E=K_M$. As Remark 3.1 (c) shows this gives an upper bound for 
$\operatorname{vol}(M)=\int_M \omega^n_\ast/n!$\,.
\hfill{\,}\qed
\enddemo

\remark{Remark 5.1}
We can prove a complete generalization of the asymptotic Morse inequalities of Demailly \cite{De1}
for the $L^2$ cohomology of the covering of a compact manifold $X$. 
For this purpose we elabotate the proof of Theorem 2.1. As there we exploit the idea of 
Witten--Demailly of constructing
a family of subcomplexes of the $L^2$--Dolbeault complex having the same cohomology.
First let us introduce cohomology.
Let us denote by $N^q(\bar{\partial})$ the kernel and 
by $R^q(\bar\partial)$ the range of
${{\bar\partial}}$ , by 
$N^q({\bar\partial}^{\ast})$ the kernel of ${\bar\partial}^{\ast}$ and by
$N^q(k^{-1}\latk)$ the kernel of $k^{-1}\latk$, all acting on $L^2_{0,q}(\tX,\widetilde E^k\otimes \tf)$
where $\tf$ is a $\g$--invariant holomorphic vector bundle of rank $r$.
We have
${\Cal H}^{0,q}_{(2)}(\tX,\widetilde E^k\otimes\tf):=N^q(k^{-1}\latk)
=N^q (\db)\cap N^q (\db^{\ast})\,,$
where the first equality is the definition of the space of harmonic forms
and the second is a consequence of the completeness of the metric.
If $q=0$ then ${\Cal H}^{0,0}_{(2)}(\tX,\widetilde E^k\otimes\tf)$ coincides
to the space of holomorphic $L^2$ sections of $\te^k\otimes\tf$.
We note also the orthogonal decomposition
$N^q(\bar\partial)=N^q(k^{-1}\latk)\oplus\overline
{R^{q-1}(\bar\partial)}$ so that 
$${\Cal H}^{0,q}_{(2)}(\tX,\widetilde E^k\otimes \tf)=N^q(\bar\partial)/\overline
{R^{q-1}(\bar\partial)}=: H^{0,q}_{(2)}(\tX,\widetilde E^k\otimes\tf)$$
the last group being the (reduced) $L^2$ cohomology. 

We apply the results of \S 1 in the following form. Since $X$ is compact we can take the 
set $\Om\Subset X$ to be $X$ so that $\widetilde\Om=\tX$. We do not use any special metric 
but take an arbitrary metric on $X$ and its pull--back on $\tX$. Moreover we have $\latk=\latd$.
Since $k^{-1}\latk$ commutes with $\db$ it follows that the spectral projections of 
$k^{-1}\latk$ commute with $\db$ too, showing thus $\db L^q_k(\lam)\subset L^{q+1}_k(\lam)$
and therefore we have a complex of $\g$--modules of finite $\g$--dimension:
$$0@>>{\;\;\;}>L^0_k(\lam)@>{\;\db_\lam\;}>>L^1_k(\lam)
@>{\;\db_\lam\;}>>\dotsm @>{\;\db_\lam\;}>>L^n_k(\lam)@>>{\;\;\;}>0\,.\eqno(5.1)$$
$k^{-1}\latk$ commutes also with $\db^\ast$ and $\left({\db_{\lam}}\right)^{\ast}$
equals the restriction of $\db^\ast$ to $ L^q_k(\lam)$. Keeping this in mind it is easy 
to see that
$$N^q(\db_\lam)/\overline {R^{q-1}(\db_\lam)}=
\Big\lbrace u\in L^q_k(\lam)\,:\,\db_{\lam} u=0\,,\;\left({\db_{\lam}}\right)^{\ast} u=0\Big
\rbrace={\Cal H}^{0,q}_{(2)}(\tX,\widetilde E^k\otimes\tf)\,.\eqno(5.2)$$
We can now apply the following lemma (see \cite{Sh}). 
\proclaim{Algebraic Lemma}
Let
$0\longrightarrow L_0 @>{\;d_0\;}>> L_1@>{\;d_1\;}>> \dotsm @>{\;d_n\;}>> L_n\longrightarrow 0$
be a complex of $\Gamma$--modules ($d_q$ commutes with the action
of $\Gamma$ and $d_{q+1}d_q=0$). If
$l_q=\dim_{\scriptscriptstyle\Gamma}L_q$ is finite and
$h_q=\dim_{\scriptscriptstyle\Gamma}H^q(L)$ where 
$H_q(L)=N(d_q)/\overline{R(d_{q-1})}$,
$$\sum\limits_{j=0}^q(-1)^{q-j}h_j\leqslant \sum\limits_{j=0}^q(-1)^{q-j}l_j
$$
for every $q=0,1,...,n$ and for $q=n$ the inequality becomes equality.
\endproclaim 

\noindent
The Algebraic Lemma for the complex (5.1) and relation (5.2) yield
$$ \sum\limits_{j=0}^{q}(-1)^{q-j}\dg H^{0,j}_{(2)}(\tX,\te^k\otimes\tf)\leqslant 
\sum\limits_{j=0}^{q}(-1)^{q-j}
N^j_{\sG}\left(\lambda,\frac{1}{k}\latk\right)$$
for $q=0,1,\dotsc,n$ and for $q=n$ the inequality becomes equality. 
We apply now (1.11):
$$\eqalign{
\sum\limits_{j=0}^{q}(-1)^{q-j}\dg H^{0,j}_{(2)}(\tX,\te^k\otimes\tf)\leqslant
k^n \big(I^q (U,\lam)&-I^{q-1} (U,\lam)+\dotsm \cr
&+(-1)^q I^0 (U,\lam)\big)+o(k^n)\,,}$$
for $k\longrightarrow\infty$. We can now let $\lam$ go to zero through values  
$\lam\in\Bbb R\setminus\Cal N$. We have thus proved the following.

\proclaim{Theorem 5.4}
Let $\tX$ be a Galois covering of group $\g$ of a compact manifold $X$.
As $k\to\infty$, the following strong Morse inequalities hold for every $q=0, 1, \dots, n$ {\rom :}
$$\sum\limits_{j=0}^{q}(-1)^{q-j}\dim_{\scriptscriptstyle\Gamma}H^{0,j}_{(2)}(\tX,
\widetilde E^k\otimes
\widetilde F)
\le r \frac{k^n}{n!} \int_{X(\leqslant q)}(-1)^q\left(\frac{\imath}{2\pi}{\boldkey c}(E)
\right)^n+o(k^n).$$
with equality for $q=n$ (asymptotic $L^2$ Riemann-Roch formula).
\endproclaim
In particular
$\dg H^0_{(2)}(\tX,\te^k\otimes\tf)\geqslant \frac{k^n}{n!}\int_{X(\leqslant 1)}
\left(\frac{\imath}{2\pi}{\boldkey c}(E)\right)^n +o(k^n)$.
It follows that if $E$ satisfies (D) then for $k\longrightarrow\infty$
$$\eqalign{\dim_{\sG}H^{0}_{(2)}(\tX,\widetilde E^k\otimes
\widetilde F)&\approx k^n\,,\cr
\dg H^q_{(2)}(\tX,\widetilde E^k\otimes F)&=o(k^n)\,,
\quad q\geqslant 1\,.}$$
Hence the usual dimension of the space of holomorphic $L^2$
sections has the same cardinal as $|\g|$ for large $k$.
This is a generalization of the result of Napier
\cite{Nap} that $\tX$ is holomorphically convex with respect to $\te^k$ for large $k$
if $X$ is projective and $E$ is positive. 
If the canonical bundle $K_X$ satisfies condition (D), i.e. if there exists a metric $\omega$
on $M$ such that $\int_{X(\leqslant 1)}(-\operatorname{Ric}\omega)^n>0$ where $X(\leqslant 1)$ 
is the set of points where $-\operatorname{Ric}\omega$ is nondegenerate and has at most one 
negative eigenvalue, then $\dg H^{0}_{(2)}(\tX,K_{\tX}^{\otimes k})\approx k^n$\,. 
\comment
This statement may be 
regarded as a converse to the theorem of Poincar\'e and Kodaira saying that if $\tX$ is a bounded open 
set in $\Bbb C^n$ on which $\g$ acts freely and properly discontinuously then $X=\tX/\g$ has ample 
canonical bundle $K_X$ (see \cite{Ko}).
\endcomment 
\endremark

\remark{Remark 5.2}
In Proposition 4.1 we have treated the case of a singular hermitian line bundle $(E,h)$ over a 
compact manifold $X$. The condition on the singularities were that they are concentrated
on an analytic set and moreover the curvature is positive near this analytic set.
Then we can work on the complement of the analytic set and by means of the basic estimate
study its coverings.
If we are interested only in the coverings of $X$ then we can rule out the condition 
of positivity near the singularities. 
Namely, when the singularities of the metric are algebraic (cf. \cite{De2}),
Bonavero \cite{Bon} shows that the Morse inequlities are true for the cohomology
of $E^k$ twisted with the corresponding sequence of Nadel's multiplier ideal sheaves.
Given a Galois covering as above we can adapt his proof to estimate the von Neuman 
dimension of the space $H^0_{(2)}(\tX,\te^k\otimes {\Cal I}_k(\widetilde h))$ of $L^2$
holomorphic sections in $\te^k$ twisted with the Nadel's multiplier ideal sheaf coming
from the singularities of the $\g$--invariant metric $\widetilde h$ on $\te^k$ 
(which is the pull--back of a Nadel multiplier ideal sheaf on $X$). 
The conclusion is that when (D) is true, the integral 
being taken over the regular set of the curvature current, then the von Neuman dimension of 
$H^0_{(2)}(\tX,\te^k\otimes {\Cal I}_k(\widetilde h))$ grows as $k^n$ for large $k$.
\endremark

\remark{Remark 5.3}
Using the approach of this section we can study the growth of the cohomology groups of coverings of
$q$--convex and $q$--concave manifolds. We can either use complete metrics or follow \cite{GHS}
and use the $\db$--Neumann problem setting. Let us give the statements in the latter set-up.
Consider a $q$--convex manifold $X$ in the sense of \cite{AG}, i.e. there exists a smooth exhausting 
function $\f:X\longrightarrow\Bbb R$ such that $\imath\partial\db\f$ has at least $n-q+1$
positive eigenvalues outside a compact set $K$ ($n=\dim X$\,, $1\leqslant q\leqslant n-1$). 
Consider $X_c=\{\f<c\}\supset K$
with smooth boundary. Then the Levi form of $bX_c$ has at least $n-q$ positive eigenvalues.
Let us consider a Galois covering $\widetilde X_d$
of a bigger sublevel set $X_d\supset X_c$ and denote
by $\widetilde X_c$ the induced covering of $X_c$. As usual we denote by $\g$ the group of 
deck transformations. Let us consider also a line bundle $E$ over $X$ and denote by $\te$
its lifting to $\widetilde X_d$. Both  $\widetilde X_c$ and $\te$ come with the liftings
of metrics defined on $X_d$.
We define the (reduced) $L^2$ cohomology groups $H^j_{(2)}(\widetilde X_c,\te^k)$ with respect to
these metrics. By \cite{GHS} we know that $\dg H^j_{(2)}(\widetilde X_c,\te^k)<\infty$ for 
$j\geqslant q$. With the method used in this paper we can prove that for $j\geqslant q$ and
$k\longrightarrow\infty$\,:
\roster
\item $\dg H^j_{(2)}(\widetilde X_c,\te^k)=O(k^n)$\,.
\item If $E$ is $q$--positive outside $K$ (its curvature has at least $n-q+1$ positive eigenvalues)
we have an explicit bound,
$$\dg H^j_{(2)}(\widetilde X_c,\te^k)\leqslant \frac{k^n}{n!}\int_{X(j)}(-1)^j
\left(\frac{\imath}{2\pi}{\boldkey c}(E)\right)^n+o(k^n).$$   
\endroster
The proof consists of showing that the basic estimate holds in bidegree $(0,j)$ on $\widetilde X_c
\subset\widetilde L$, where $L$ is a compact set of $X_c$\,, for forms satisfying the $\db$--Neumann 
conditions on $b\widetilde X_c$\,. This is achieved using the liftings of the metrics constructed in \cite{AV}
where the case $\g$ trivial is treated. Then we can apply again the analysis from \S 1. 
If $E$ is $q$--positive outside $K$ then the leading term in (1) simplifies as shown in (2). 
These estimates were obtained in the case $\g=\{\Id\}$ in 
\cite{Bou} for certain complete metrics on $X_c$ which permit to prove the same inequalities for
the full cohomology group $H^j(X_c,E^k)$. For the case of coverings we have to restrict ourselves to $L^2$
cohomology groups.  As for coverings of $q$--concave manifolds we get the same conclusion as in (1) for 
$j\leqslant n-q-1$. The nice simplification of the leading term holds if we impose a negativity condition
outside a compact set. However there are cases of concave manifolds and positive bundles for which  
we have an effective estimate of $\dg H^0_{(2)}(\widetilde X_c,\te^k)$, see \S 3. 
\endremark 

\subhead
\S 6 Weak Lefschetz theorems
\endsubhead

Nori \cite{No} generalized the Lefschetz hypersurface theorem. Assume $X$ and $Y$ are smooth connected
projective manifolds and $Y$ is a hypersurface in $X$ with positive normal bundle and $\dim Y\geqslant 1$.
Then the image of $\pi_1(Y)$ in $\pi_1(X)$ is of finite index. Recently, Napier and Ramachandran \cite{NR}
proposed an analytic approach and generalized Nori's theorem showing that $Y$ may have arbitrary codimension
(but $\dim Y\geqslant 1$). They use the $\db$--method on complete K\"ahler manifolds to separate 
the sheets of appropriate coverings. In the sequel we use the Riemann--Roch inequalities to study
non--necessarily K\"ahler manifolds. 
However our method requires that the image group is normal since we can deal only with 
Galois coverings.
First we introduce the notion of formal completion.
Let $Y$ be a complex analytic subspace of the manifold $U$ and denote by ${\Cal I}_Y$
the ideal sheaf of $Y$. The formal completion $\widehat U$ of $U$ with respect to $Y$ is the 
ringed space $({\widehat U},{\Cal O}_{\widehat U})=(Y,\projlim {\Cal O}_{U}/ {\Cal I}^{\nu}_Y)$.
If $\Cal F$ is an analytic sheaf on $U$ we denote by $\widehat {\Cal F}$ the sheaf  
$\widehat {\Cal F}=\projlim {\Cal F}\otimes ({\Cal O}/{\Cal I}^{\nu}_Y)$. If ${\Cal F}$ is coherent then
$\widehat {\Cal F}$ is too. Moreover by Proposition VI.2.7 of \cite{BS} the kernel of
the mapping $H^0(U,{\Cal F})\longrightarrow H^0(\widehat U,\widehat {\Cal F})$ consists of the sections of 
${\Cal F}$ which vanish on a neighbourhood of $Y$. Hence for locally free $\Cal F$ the map is 
injective.
\proclaim{Theorem 6.1}
Let $M$ be a hyper $1$--concave manifold carrying a line bundle $E$ which satisfies {\rm (D)} 
and is semi-positive
outside a compact set. Let $Y$ be a connected compact complex subspace of $M$  satisfying: 
{\rm (i)} for any $k$, 
$\dim H^0 (\widehat M,{\widehat {\Cal F}_k})
<\infty$, where ${\Cal F}_k={\Cal O}(E^k\otimes K_M)$, {\rm (ii)} the image $G$ of 
$\pi_1(Y)$ in $\pi_1(X)$ is normal in $\pi_1(X)$. 
Then $G$ is of finite index in $\pi_1(X)$.
\endproclaim\par
\demo{Proof}
We follow the proof given in \cite{NR}. Since $G$ is normal there exists a connected 
Galois covering $\pi:\tm\longrightarrow M$ such that the group of deck transformations is $\g=\pi_1(M)/G$. 
The cardinal $|\g|$ equals the index of $G$ in $\pi_1(M)$. Let $\te=\pi^{-1}E$. By Theorem 3.1,
there exists $C>0$ such that for large $k$,
$\dg H_{(2)}^{n,0}(\tm,\te^k)\geqslant C\,k^n$. Let us choose a small open neighbourhood
$V$ of $Y$ such that $\pi_1(Y)\longrightarrow\pi_1(V)$ is an isomorphism; so the image of $\pi_1(V)$ in $\pi_1(M)$
is $G$. Hence, if we denote by $\jmath$ the inclusion
of $V$ in $M$, there exists a holomorphic lifting $\widetilde \jmath:V\longrightarrow\tm$, 
$\pi\circ\widetilde \jmath=\jmath$. Since $\widetilde\jmath$ is locally biholomorphic the pull--back map
${\widetilde\jmath}^{\,*}: H_{(2)}^{n,0}(\tm,\te^k)\longrightarrow H^{n,0}(V,E^k)$
is injective. On the other hand $H^0(V,{\Cal F}_k)\hookrightarrow H^0(\widehat V,\widehat{\Cal F}_k)
=H^0(\widehat M,\widehat{\Cal F}_k)$. By (i) the latter space is finite dimensional so 
$\dim H^{n,0}_{(2)}(\tm,\te^k)<\infty$. We know that
$\dg H_{(2)}^0(\tm,\te^k\otimes K_{\tm})>0$ for $k>C^{-1/n}$. If $\g$ were infinite this would yield 
$\dim H^{n,0}_{(2)}(\tm,\te^k)=\infty$ which is a contradiction. Therefore $|\g|<\infty$ and
$\dim H^{n,0}_{(2)}(\tm,\te^k)\geqslant C\,|\g|\,k^n\geqslant |\g|$ for $k>C^{-1/n}$.
Thus $|\g|\leqslant \dim H^0(\widehat M,\widehat{\Cal F}_k)$ for large $k$. 
\hfill{\,}\qed\enddemo\par
\remark{Remark 6.2}

\noindent
(a) By a theorem of Grothendieck \cite{Gro}, condition (i) is fulfilled if $Y$ is locally a complete intersection
with ample normal bundle $N_Y$ (or $k$--ample in the sense of Sommese, $k=\dim Y-1$).

\noindent
(b) We can replace condition (i) with the requirement that $Y$ has a fundamental system of 
pseudoconcave neighbourhoods $\{V\}$. Then $\dim H^0 (V,{\Cal F}_k)$ is finite by \cite {An}.
This happens for example if $Y$ is a smooth hypersurface and $N_Y$ has at least
one positive eigenvalue or, if $Y$ has arbitrary codimension, if $N_Y$ is sufficiently 
positive in the sense of Griffiths \cite{Gri1}.

\noindent
(c) Condition (ii) is trivially satisfied if $\pi_1(Y)=0$. Thus, if $M$ contains a
simply connected subvariety satisfying either (a) or (b), $\pi_1(M)$ is finite.

\noindent
(d) By Corollary 3.6, Theorem 6.1 can also be applied to the pertubed structures considered there.
\endremark

Using Proposition 4.3 we can can show that Nori's theorem holds for
all Moishezon spaces $X$. 
\proclaim{Theorem 6.2}
Let $X$ be an irreducible reduced normal Moishezon compact
complex space and let $E$ be the (positive in the sense of currents) line
bundle given by Theorem {\rm 4.2}. Suppose that $M$ is a Zariski open
set of $X$ and $Y\subset\operatorname{Reg}(M)$ be a connected
compact complex subspace such that for any $k$, 
$\dim H^0 (\widehat M,{\widehat {\Cal E}_k})<\infty$, where ${\Cal E}_k={\Cal O}(E^k)$.
Then the image $G$ of $\pi_1(Y)$ in $\pi_1(M)$ is of finite index in $\pi_1(M)$. 
\endproclaim
\demo{Proof}
Since $X$ is normal we have an isomorphism
$\pi_1(\operatorname{Reg}M)\longrightarrow\pi_1(M)$, so that we may
assume $M\subset\rg$.
We find a connected unramified covering $p:\tm\longrightarrow M$ such
that $p_{\ast}\pi_1(\tm)=G$. If $d$ is the number of sheets,
$d=|\pi_1(M)/G|$, the index of $G$ in $\pi_1(M)$. The preceding proof 
applies by using Proposition 4.3 instead of Theorem 3.1. and the usual
dimension instead of the $\g$--dimension.
\hfill{\,}\qed
\enddemo
Note that Napier and Ramachandran also
considered cases when $X$ is not necessarily projective, but their
result does not imply diectly Theorem 6.2.
 

\enddocument

\Refs
\widestnumber\key{RRRR}

\ref
\key An
\by A. Andreotti
\paper Th\'eor\`emes de d\'ependance alg\'ebrique sur les espaces complexes
pseudo-concaves
\jour  Bull.\ Soc.\ Math.\ France 
\vol 91
\yr 1963
\pages 1--38
\endref

\ref
\key At
\by M. F. Atiyah
\paper Elliptic operators, discrete groups and von Neumann
algebras
\inbook Ast\'erisque
\yr 1976
\vol 32--33
\pages 43--72
\endref

\ref
\key AG
\by A. Andreotti, H. Grauert
\paper Th\'eor\`eme de finitude pour la
cohomologie des espaces complexes  
\jour Bull.\ Soc.\ Math.\ France 
\vol 90
\yr 1962
\pages 193--259 
\endref

\ref
\key AV
\by A. Andreotti, E. Vesentini 
\paper Carleman estimates for the
Laplace-Beltrami equation on
complex manifolds
\jour Inst.\ Hautes\ Etudes\ Sci.\ Publ.\ Math.\ No.
\vol 25
\yr 1965 
\pages 81--130
\endref

\ref
\key AT\by A. Andreotti, G. Tomassini
\paper Some remarks on pseudoconcave manifolds
\inbook Essays in Topology and Related Topics, dedicated to G. de Rham
\publ Springer \yr 1970 \pages 84--105
\publaddr Berlin--Heidelberg--New York
\eds R. Narasimhan, A. Haefliger 
\endref

\ref
\key Bou
\by T. Bouche
\paper Inegalit\'es de Morse pour la $d''$--cohomologie sur
une vari\'et\'e
non--compacte
\jour Ann.\ Sci.\ Ecole\ Norm.\ Sup. 
\vol 22 
\yr 1989 
\pages 501--513
\endref

\ref
\key BS
\by C. B\u{a}nic\u{a}, O. St\u{a}n\u{a}\c{s}il\u{a}
\book Algebraic methods in the global theory of complex spaces
\publ Wiley\publaddr New York\yr 1976
\endref

\ref B
\key B
\by N. Bourbaki
\book Vari\'et\'es diff\'erentielles et analytiques
\moreref Fascicule de r\'esultats
\yr 1967
\publ Hermann\publaddr Paris
\bookinfo Actualit\'es scientifiques et industrielles 1333
\endref


\ref
\key CFKS
\by H.L. Cycon, R.G. Froese, W. Kirsch, B. Simon
\book Schr\"odinger operators with applications to quantum physics
\publ Springer--Verlag\yr 1987
\bookinfo Text and Monographs in Physics  
\endref


\ref
\key De1
\by J. P. Demailly
\paper Champs magn\'etiques et inegalit\'es de Morse pour la
$d''$--cohomologie
\jour Ann.\ Inst.\ Fourier 
\vol 35 
\yr 1985
\pages 189--229
\endref
 
\ref
\key De2
\by J. P. Demailly
\paper Regularization of closed positive currents and intersection theory
\jour J.\ Alg.\ Geom.
\vol 1\yr 1992\pages 361--409
\endref

\ref
\key GR
\by H. Grauert, O. Riemenschneider
\paper Verschwindungss\"atze f\"ur analytische
Kohomologiegruppen auf komplexen R\"aumen
\jour Invent. Math. 
\vol 11
\yr 1970
\pages 263--292
\endref

\ref 
\key Gri1
\by Ph. A. Griffiths
\paper The extension problem in complex analysis; embedding with positive 
normal bundle
\jour Amer. J. Math. \vol 88 \yr 1966 \pages 366--446
\endref

\ref 
\key Gri2
\by Ph. A. Griffiths
\paper Complex--analytic properties of certain Zariski open sets on algebraic varieties
\jour Ann. of Math.\vol 94\yr 1971 \pages 21--51
\endref

\ref
\key Gro
\by A. Grothendieck
\book Cohomologie locale des faisceaux coh\'erents et th\'eor\`emes
de Lefschetz locaux et globaux
\publ North--Holland\publaddr Amsterdam\yr 1968
\endref

 
\ref
\key GHS
\by M. Gromov, M. G. Henkin, M. Shubin
\paper $L^2$ holomorphic functions on pseudo--convex coverings
\jour Contemp. Math. \vol 212\yr 1998 \pages 81--94
\publaddr Providence, RI
\publ AMS 
\endref

\ref
\key He
\by G.Henniart
\paper Les in\'egalit\'es de Morse (d'apr\`es Witten)
\jour Ast\'erisque\issue 121/122\vol 1983/84
\pages 43--61\yr 1985
\endref

\ref
\key Ko
\by J. Koll\'ar
\book Shafarevich maps and automorphic forms
\publ Princeton University Press
\publaddr Princeton, NJ
\yr 1995
\endref

\ref
\key Le
\by L. Lempert
\paper Embeddings of three dimensional Cauchy--Riemann manifolds
\jour Math. Ann.\vol 300\pages 1-15\yr 1994
\endref

\ref
\key Ma
\by G. Marinescu
\paper Asymptotic Morse Inequalities for Pseudoconcave Manifolds
\jour Ann.\ Scuola\
Norm.\ Sup.\ Pisa\ Cl.\ Sci.\ (4) 
\vol 23 
\yr 1996 
\issue 1
\pages 27--55
\endref  

\ref
\key NT
\by A. Nadel, H. Tsuji
\paper Compactification of complete K\"ahler manifolds
of negative Ricci curvature\jour J. Differential Geom. \vol 28 \yr 1988
\issue 3\pages  503--512\endref 

\ref
\key Nad
\by A. Nadel
\paper On complex manifolds which can be compactified by
adding finitely many points
\jour Invent. Math. \vol 101 \yr 1990 \issue 1
\pages 173--189 
\endref

\ref
\key Nap
\by T. Napier \paper Convexity properties of coverings of smooth projective
varieties\jour  Math. Ann. \vol 286 \yr 1990 \pages 433--479.
\endref

\ref
\key NR
\by T. Napier, M. Ramachandran
\paper The $L^2$--method, weak Lefschetz theorems and the topology of K\"ahler
manifolds
\jour JAMS\vol 11\issue 2\pages 375--396
\endref 

\ref
\key No
\by M.V. Nori\paper Zariski's conjecture and related problems
\jour Ann. Sci. Ec. Norm.Sup.\vol 16\yr 1983\pages 305--344
\endref

\ref
\key Oh1
\by T. Ohsawa\paper Hodge spectral sequence and symmetry on compact K\"ahler
spaces\jour Publ. Res. Inst. Math. Sci. \vol 23\yr 1987
\pages 613--625
\endref

\ref
\key Oh2
\by T. Ohsawa\paper Isomorphism theorems for cohomology groups of weakly
$1$--complete manifolds
\jour Publ. Res. Inst. Math. Sci. \vol 18\yr 1982
\pages 191--232
\endref
 
\ref
\key RSS
\by G.V. Rozenblum, M.A. Shubin, M.Z. Solomyak
\paper Spectral Theory of Differential operators
\inbook Partial Differential Equations VII, Encyclopedia of Mathematical Sciences,
\vol 64\yr 1994 \publ Springer--Verlag
\endref

\ref
\key Sa
\by L. Saper
\paper $L^2$-cohomology and intersection homology of 
certain algebraic
varieties with isolated singularities
\jour Invent. Math. \vol 82 \yr 1985
\issue 2\pages 207--255
\endref

\ref
\key Si1
\by Y.~T. Siu 
\paper A vanishing theorem for semipositive line bundles over non-K\"ahler manifolds
\jour J. Diff. Geom.\vol 19
\pages 431--452 
\yr 1984
\endref 

\ref
\key Si2
\by Y.~T. Siu 
\paper Some recent results in complex manifold theory related
to vanishing theorems for the semipositive case
\inbook Lecture Notes in Math.
\pages 169--192 
\vol 1111
\publ Springer
\publaddr Berlin-New York
\yr 1985
\bookinfo Workshop Bonn 1984
\endref 

\ref
\key SY\by Y. T. Siu, S. T. Yau
\paper Compactification of negatively curved complete K\"ahler manifolds of 
finite volume
\jour Ann. Math. Stud.\vol 102\pages 363--380\yr 1982
\endref

\comment
\ref\key Si2
\by Y. T. Siu
\paper A vanishing theorem for semipositive line bundles over non--K\"ahler
manifolds
\jour J. Differential Geometry
\vol 20 \yr 1984 \pages 431--452
\endref
\endcomment

\ref
\key Sh
\by M. Shubin
\paper Semiclassical asymptotics on covering manifolds and Morse
inequalities
\jour Geom. Funct. Anal.
\vol 6
\yr 1996
\issue 2
\pages 370--409
\endref

\ref
\key Ta
\by S. Takayama
\paper A differential geometric property of big line
bundles
\jour T\^ohoku Math. J.\ (2) \vol 46 
\yr 1994 \issue 2\pages 281--291
\endref

\comment
\ref
\key TC
\by R. Todor, I. Chiose
\paper Holomorphic Morse inequalities on covering manifolds
\jour Preprint xxx.lanl.gov/abs/math.CV/9802110 
\endref
\endcomment

\ref
\key Zu\by S. Zucker\paper Hodge theory with degenerating coefficients: 
$L^2$ cohomology in the Poincar\'e metric
\jour Ann. Math.\vol 109\pages 415--476\yr 1979
\endref

\ref
\key Wi
\by E. Witten
\paper Supersymmetry and Morse theory
\jour J. Diff. Geom.\vol 17\pages 661--692
\yr 1982
\endref



\endRefs
\bye